\documentclass[12,a4paper]{amsart}
\usepackage{amscd,amssymb,array,amsmath,epsfig,booktabs}
\usepackage{graphicx,color}

\def\limsup{\mathop{\overline{\rm lim}}}

\newtheoremstyle{sty}%
{}%
{}%
{\it}%
{}%
{\bf}%
{}%
{.5em}%
{}%
\theoremstyle{sty}
\newcommand{\D}{\displaystyle}
\newtheorem{theorem}{Theorem}[section]
\newtheorem*{theoremA}{Theorem A}
\newtheorem{corollary}[theorem]{Corollary}
\newtheorem{proposition}[theorem]{Proposition}
\newtheorem{remark}[theorem]{Remark}
\newtheorem{lemma}[theorem]{Lemma}
\newtheorem{example}[theorem]{Example}

\numberwithin{equation}{section}

\renewcommand{\Re}{{\rm Re}\,}
\renewcommand{\Im}{{\rm Im}\,}
\newtheorem*{note}{Note added in Proof of Theorem 3.4}
\allowdisplaybreaks[4]

\begin{document}

\begin{center}
	\huge{\textbf{Subnormal solutions of non-homogeneous periodic ODEs, Special functions  and related  polynomials}}$^\dag$
\end{center}
\title{}
\date{$21^{\textrm{st}}$ June 2010}

\author[Y. M. Chiang]{Yik-Man Chiang}

\address{Department of Mathematics, The Hong Kong
University of Science \& Technology, Clear Water Bay, Kowloon, Hong
Kong, China} \email{machiang@ust.hk}

\author[K. W. Yu]{Kit-Wing Yu}
\address{Department of Mathematics, The Hong Kong
University of Science \& Technology, Clear Water Bay, Kowloon,
Hong Kong, China} \email{kitwing@ust.hk}

\begin{abstract}
This paper offers a new and complete description of subnormal
solutions of certain non-homogeneous second order periodic
linear differential equations first studied by Gundersen and
Steinbart in 1994. We have established a previously unknown
relation that the general solutions (\textit{i.e.}, whether
subnormal or not) of the DEs can be solved explicitly in terms of classical special functions, namely the Bessel, Lommel and Struve functions, which are important because of their numerous physical applications. In particular, we
show that the subnormal solutions are written explicitly in
terms of the degenerate Lommel functions $S_{\mu,\,\nu}(\zeta)$ and several classical special polynomials related to the Bessel functions. In fact, we solve an equivalent problem in special functions that each branch of the Lommel function $S_{\mu,\, \nu}(\zeta)$ degenerates if and only if $S_{\mu,\, \nu}({\rm e}^z)$ has finite order of growth in $\mathbf{C}$. We achieve this goal
by proving new properties and identities for these functions.  A number of semi-classical quantization-type results are obtained as consequences. Thus our
results not only recover and extend the result of Gundersen and Steinbart \cite{GS94}, but the new
identities and properties found for the Lommel functions are of
independent interest in a wider context.
\end{abstract}


\dedicatory{Dedicated to the sixtieth birthday of Pit-Mann Wong.}


\thanks{$^\dag$ This research was partially supported by the
University Grants Council of the Hong Kong Special Administrative
Region, China (HKUST6145/01P and 600806)}

\maketitle

\tableofcontents

\section{\bf Introduction and main results}\label{S:introduction}
The growth of entire solutions of the equation

  \begin{equation}\label{E:Gundersen-Steinbart}
    f^{\prime\prime}+P({\rm e}^z)f^{\prime}+Q({\rm e}^z)f=R_1({\rm e}^z)+R_2({\rm e}^{-z}),
  \end{equation}\smallskip

\noindent where $P(\zeta),\, Q(\zeta),\, R_1(\zeta)$
and $R_2(\zeta)$ are polynomials in $\zeta$ and that $P(\zeta)$ and
$Q(\zeta)$ are not both constant, were considered in \cite{Ga91}, \cite{GS92},
\cite{GS94}, \cite{GSW98} and \cite{La93}. It was shown that certain
\textit{subnormal solutions} can be written in the form $f(z)={\rm
e}^{dz}S({\rm e}^z)$ where $S(\zeta)$ is a polynomial and $d$ is a constant. The same problem when the equation (\ref{E:Gundersen-Steinbart}) is homogeneous was considered, for examples, in \cite{BS97}, \cite{BL83}, \cite{BLL86} and \cite{Fr61}.\smallskip

In this paper, we exhibit a previously unknown relation that the
solutions of a subcase of the equation
(\ref{E:Gundersen-Steinbart}) when $\deg P<\deg Q\le 1$ can be solved in terms of the sum of
the Bessel functions and the Lommel function
$S_{\mu,\,\nu}(\zeta)$. In the most general consideration the existence of the
subnormal solutions of this important subclass of
(\ref{E:Gundersen-Steinbart}) is equivalent to the degeneration of the $S_{\mu,\,\nu}(\zeta)/\zeta^{\mu-1}$ into a polynomial in $\zeta$ and $1/\zeta$. In several specialized considerations, classical special polynomials related to the Bessel functions such as the Struve
functions, the Neumann polynomials, the Gegenbauer polynomials and
the Schl\"afli polynomials \cite{Wa44}, \S 9.1-9.3 and \S 10.4  are needed in order to describe the subnormal solutions.
\smallskip

The Lommel functions $S_{\mu,\,\nu}(\zeta)$ have numerous applications in, for examples,
electromagnetic scattering in a multilayered medium \cite{CY2000}, \cite{CY2001},
thermal inflation \cite{LS1996}, one-dimensional stochastic model with branching and coagulation reactions
\cite{MB2001}, oscillatory limited compressible fluid flow \cite{SGD2003},
computation of toroidal shells and propeller blades \cite{Tumarkin1959} and strain gradient elasticity theory for antiplane shear cracks \cite{UA2000}, etc. In particular, its special case, when
$\mu=\nu$, the Struve function ${\bf H}_\nu(\zeta)$ also occurs in
many applications \cite{Wa44}, pp. 328--338. See for examples, in the theory of loud speakers \cite{McLachlan1934} and in the theory of light \cite{Walker1904}, { chap. 7}.
\smallskip

An entire solution $f(z)$ of (\ref{E:Gundersen-Steinbart}) is called {\it
subnormal} if either

    \begin{equation}\label{E:subnormal}
    \limsup_{r \to+\infty}{\log\log M(r,\, f)\over
    r}=0\quad\textrm{or}\quad \D \limsup_{r \to+\infty}\frac{\log T(r,\,
    f)}{r}=0
    \end{equation}\smallskip

\noindent holds. Here $\D M(r,\, f)=\max_{|z|\le r}|f(z)|$ denotes
the usual maximum modulus of the entire function $f(z)$ and $T(r,\,
f)$ is the Nevanlinna characteristics of $f(z)$. We denote the
the {\it order} of a meromorphic function $f(z)$ by
$\sigma(f)=\limsup_{r\to+\infty}\log\log M(r,\, f)/\log r=\limsup_{r\to+\infty}\log T(r,\, f)/\log r$. We refer the reader to \cite{Ha75} or \cite{La93} for
the details.\smallskip

In \cite{GS94}, Gundersen and Steinbart proved

    \begin{theoremA}\label{T:Gundersen-Steinbart}
    Suppose that $\deg P <\deg Q$ in the non-homogeneous differential
    equation {\rm (\ref{E:Gundersen-Steinbart})}
    and that {\rm (\ref{E:Gundersen-Steinbart})} admits a
    subnormal solution $f(z)$. Then $f(z)$ must have the form

        \begin{equation}
        \label{E:Gundersen-Steinbart-soln} \D f(z)=S_1({\rm e}^z)+S_2({\rm e}^{-z}),
        \end{equation}\smallskip

    \noindent where $S_1(\zeta)$ and $S_2(\zeta)$ are polynomials in $\zeta$.
    \end{theoremA}

Gundersen and Steinbart also considered the cases when $\deg P >\deg Q$ and  $\deg P = \deg Q$, respectively, and obtained subnormal solutions similar to (\ref{E:Gundersen-Steinbart-soln}) under the same assumption that the equation (\ref{E:Gundersen-Steinbart}) admits a subnormal solution. In this paper we shall only consider the case $\deg P<\deg Q$ and we refer the reader to \cite{GS94} for other details.
\smallskip

\begin{remark} \rm We note that Theorem A and the other results obtained by Gundersen and Steinbart \cite{GS94} are generalizations to those of Wittich's \cite{Wi67} for the periodic
homogeneous equation

    \begin{equation}\label{E:Wittich}
    f^{\prime\prime}+P({\rm e}^z)f^{\prime}+Q({\rm e}^z)f=0.
    \end{equation}\smallskip

\noindent Wittich showed that each subnormal solution $f(z)$ to this equation admits a representation of the form $f(z)={\rm e}^{dz}S({\rm e}^z)$, where $S(\zeta)$ is a polynomial in $\zeta$ and $d$ is a constant.
\end{remark}
\smallskip

Ismail and one of the authors showed in \cite{CI02(2)}, Remark 1.11 (see also \cite{CI02(1)}) that the subnormal
solutions of a subclass of \textit{homogeneous} differential equations of
(\ref{E:Wittich}), first considered by Frei in \cite{Fr61} and
then by Bank and Laine in \cite{BL83}, was in fact an important diatomic molecule model in Wave (quantum) mechanics proposed by P. M. Morse {in a landmark paper} \cite{Morse1929}  in 1929, and it can be solved explicitly in terms of a class of confluent hypergeometric functions -- the Coulomb Wave functions and a class of orthogonal polynomials --
the Bessel polynomials \cite{CI02(2)}. The study appears to be the first of its kind concerning subnormal solutions of special cases of (\ref{E:Wittich}) and special functions. We continue the study in this paper and to show
that special functions solutions also exist for a subclass of (\ref{E:Gundersen-Steinbart}) when $\deg P <\deg Q$. We are able to {\it solve} the equation (\ref{E:Gundersen-Steinbart}) by finding explicit solutions in terms of several classes of classical special functions, whether the solution is subnormal or not, and from which necessary and sufficient conditions for the existence of subnormal solutions can be derived as special cases
of our main results. These characterization results can be considered as a semi-classical quantization-type results for non-homogeneous equations. This matter will be further discussed in \S\ref{S:1/16-results}. \smallskip

We shall first state our main results on subnormal solutions, followed by the main result on the growth of $S_{\mu,\,\nu}({\rm e}^z)$ which is the crux of the paper.

    \begin{theorem}\label{T:Chiang-Yu-special}
    Let $\sigma,\, \nu$ be arbitrary complex constants with $\sigma$ non-zero, and let $f(z)$ be an entire solution to the differential equation

        \begin{equation}\label{E:Chiang-Yu-special}
        f''+({\rm e}^{2z}-\nu^2)f=\sigma {\rm e}^{(\mu+1) z}.
        \end{equation}\smallskip

    \begin{enumerate}
    \item[(a)] Then the general solution $f(z)$ to {\rm (\ref{E:Chiang-Yu-special})} is given by

        \begin{equation}\label{E:Chiang-Yu-special-soln-1}
        f(z)=AJ_\nu({\rm e}^z)+BY_\nu({\rm e}^z)+\sigma S_{\mu, \,\nu}({\rm e}^z),
        \end{equation}\smallskip

        \noindent where $A$ and $B$ are constants.\smallskip

    \item[(b)] The function $f(z)$ given in {\rm (\ref{E:Chiang-Yu-special-soln-1})} is
    subnormal if and only if $A=B=0$, and either

    \begin{equation*}
    \mu + \nu=2p+1\quad \mbox{or}\quad \mu-\nu=2p+1
    \end{equation*} \smallskip

    \noindent holds for a non-negative integer $p$ and      \begin{equation}\label{E:special-lommel-sol}
      S_{\mu,\, \nu}(\zeta)=\zeta^{\mu-1}\Bigg[\sum_{k=0}^p \frac{(-1)^k
      c_k}{\zeta^{2k}}\Bigg],
      \end{equation}\smallskip

    \noindent where the coefficients $c_k,\, k=0,\,1,\ldots,\, p,$ are defined by
\begin{equation}
\label{ck}
    c_0=1,\quad  c_k=\prod_{m=1}^k [(\mu-2m+1)^2-\nu^2].
\end{equation}
    \end{enumerate}
    \end{theorem}
\medskip

The \textit{Lommel function}  $S_{\mu,\, \nu}(\zeta)$ (see \S 3.1)  with respect to parameters $\mu,\, \nu$  appears to be first studied by Lommel \cite{Lo1876} in 1876. It is a particular solution of the
non-homogeneous Bessel differential equation

\begin{equation} \label{E:lommel}
    \zeta^2y^{\prime\prime}(\zeta)+\zeta y^\prime(\zeta)
    +(\zeta^2-\nu^2)y(\zeta)=\zeta^{\mu+1}.
\end{equation}\smallskip

\noindent The functions $J_\nu(\zeta)$ and $Y_\nu(\zeta)$ in
(\ref{E:Chiang-Yu-special-soln-1}) are the standard \textit{Bessel
functions of the first and second kinds} respectively. They are
two linearly independent solutions of the corresponding
homogeneous differential equation of (\ref{E:lommel}).
\smallskip

When $\mu=\nu$, we have the following special case:

    \begin{corollary} \label{T:Chiang-Yu-Struve}
    Let $f(z)$ be an entire solution to the differential equation

        \begin{equation}\label{E:Chiang-Yu-Struve}
        f''+({\rm e}^{2z}-\nu^2)f=\sigma {\rm e}^{(\nu+1) z}=\frac{{\rm e}^{(\nu+1) z}}{2^{\nu-1}\sqrt{\pi}\Gamma(\nu+\frac{1}{2})}.
        \end{equation}\smallskip

    \begin{enumerate}
    \item[(a)] Then the general solution $f(z)$ to {\rm (\ref{E:Chiang-Yu-Struve})} is given by

        \begin{equation}\label{E:Chiang-Yu-Struve-soln-1}
        f(z)=AJ_\nu({\rm e}^z)+BY_\nu({\rm e}^z)+ {\bf K}_{\nu}({\rm e}^z),
        \end{equation}\smallskip

        \noindent where $A$ and $B$ are constants.

    \item[(b)] The function $f(z)$ given in {\rm (\ref{E:Chiang-Yu-Struve-soln-1})} is subnormal if
    and only if $A=B=0$, $\nu=p+\frac12$ for a  non-negative integer $p$  \emph{and}

        \begin{equation*}\label{E:Struve-sol}
        {\bf K}_{p+\frac{1}{2}}(\zeta)=\frac{\zeta^{p-\frac{1}{2}}}{2^{p-\frac{1}{2}}\sqrt{\pi}\Gamma(p+1)}
        \cdot\Bigg[\sum_{k=0}^{p} \frac{(-1)^k
        c_k}{\zeta^{2k}}\Bigg],
        \end{equation*}\smallskip

    \noindent where each of the coefficient $c_k$ is defined in {\rm (\ref{ck})}.
    \end{enumerate}
    \end{corollary}

The function $\mathbf{K}_\nu(\zeta)$ that appears in (\ref{E:Chiang-Yu-Struve-soln-1}) above is related to the
\textit{Struve function of order} $\nu$, $\mathbf{
H}_\nu(\zeta)$, which is a particular solution of the differential
equation (\ref{E:struve-eqn}). Detailed relations
amongst ${\bf K}_\nu(\zeta)$, ${\bf H}_\nu(\zeta)$ and
$S_{\nu,\,\nu}(\zeta)$ will be given in \S \ref{Adx:Special}.5. The Struve function
was first studied by Struve in 1882 (see \cite{Wa44}, \S10.4). We shall derive analytic continuation formulae for the Lommel functions. Then the corresponding formula for the Struve function follows
as a corollary which is stated in \S B.5. Besides, we also
prove that $\mathbf{H}_\nu({\rm e}^z)$ is of infinite order of
growth (and in fact, not subnormal) for any choice of $\nu$ as a corollary of Proposition 4.1.
\medskip

We shall obtain the Theorem 1.2 as a special case of the following
more general result. We first introduce a set of more general
coefficients.
\medskip

Suppose that $n$ is a positive integer and $A,\, B,\, \nu,\, L,\,
M,\, N,\, \sigma,\, \sigma_i,\, \mu_j$ are complex numbers such
that $L,\, M$ are non-zero and at least one of  $\sigma_j$,
$j=1,\, 2,\ldots,\,n$, being non-zero. We also let

\begin{equation}
    c_{0,\,j}=1,\quad c_{k,\, j}=\prod_{m=1}^k
[(\mu_j-2m+1)^2-\nu^2]\label{ckj}
\end{equation}
\medskip
for $j=1,\, 2,\,\ldots,\, n$.

    \begin{theorem}\label{T:Chiang-Yu} Let $f(z)$ be an entire solution to the
    differential equation

        \begin{align} \label{E:Chiang-Yu}
        f''+2Nf'+\big[L^2M^2{\rm e}^{2Mz} &+(N^2 -\nu^2M^2)\big]f \notag\\
        &=\sum_{j=1}^n \sigma_j L^{\mu_j+1}M^2{\rm e}^{[M(\mu_j+1)-N]z}.
        \end{align}\smallskip

    \begin{enumerate}
      \item[(a)] Then the general solution $f(z)$ to {\rm (\ref{E:Chiang-Yu})} is given by

        \begin{equation}\label{E:Chiang-Yu-soln-1}
          f(z)={\rm e}^{-N z}\left[AJ_\nu(L {\rm e}^{M z})+ B Y_\nu (L {\rm e}^{M
        z})+\sum_{j=1}^n \sigma_j S_{\mu_j,\, \nu} (L {\rm e}^{M z})\right].
        \end{equation}\smallskip

      \item[(b)] If all the $\Re(\mu_j)$ are distinct, then the function $f(z)$ given in
      {\rm (\ref{E:Chiang-Yu-soln-1})} is  subnormal if and only if $A=B=0$ and for each non-zero $\sigma_j$,
      we have either

      \begin{equation}\label{E:conditions}
      \mu_j + \nu=2p_j+1 \quad \mbox{or}\quad \mu_j-\nu=2p_j+1,
      \end{equation}\smallskip

      \noindent where $p_j$ is a non-negative integer
      and

        \begin{equation}\label{E:lommel-sol}
        S_{\mu_j,\, \nu}(\zeta)=\zeta^{\mu_j-1}\left[\sum_{k=0}^{p_j}
        \frac{(-1)^k c_{k,\, j}}{\zeta^{2k}}\right],
        \end{equation}\smallskip

        \noindent for $j \in \{1,\,2,\ldots,\,n\}$, and where each coefficient
        $c_{k,\,j}$ is defined in {\rm (\ref{ckj})}.
    \end{enumerate}
    \end{theorem}

\begin{remark}\label{R:bessel-different-branches} \rm We  note that for all values of $\mu_j$ and $\nu$, $J_\nu(L{\rm
e}^{Mz}),\,Y_\nu(L{\rm e}^{Mz})$ and $S_{\mu_j,\,\nu}(L{\rm
e}^{Mz})$ are entire functions in the complex
$z$-plane. Hence it is a single-valued function and so is \textit{independent of the branches} of $S_{\mu_j,\,\nu}(\zeta)$.
\end{remark}

\begin{remark}\label{R:relation-theoremA-ChiangYu} \rm If we choose $L=2, M=\frac{1}{2}$ in Theorem \ref{T:Chiang-Yu} and let

\begin{align*}
P(\zeta)=2N,\, Q(\zeta)=\zeta+\bigg(N^2-\frac{\nu^2}{4}\bigg),\,
R(\zeta)=\sum_{j=1}^n
\sigma_j2^{\mu_j-1}\zeta^{\frac{1}{2}(\mu_j+1)-N},
\end{align*}\smallskip

\noindent then $P(\zeta)$ and $Q(\zeta)$ are polynomials in $\zeta$ in (\ref{E:Gundersen-Steinbart})
with $\deg P<\deg Q$. But we note that $R(\zeta)$ so chosen may not
necessary be a polynomial in $\zeta$. It follows from Theorem
\ref{T:Chiang-Yu} that any subnormal solution $f(z)$ given by
(\ref{E:Chiang-Yu-soln-1}) has the form
(\ref{E:Gundersen-Steinbart-soln}). Thus our results
(\ref{E:Chiang-Yu-soln-1}) and (\ref{E:lommel-sol}) generalize the Theorem A and give explicit formulae of the
subnormal solutions of (\ref{E:Gundersen-Steinbart-soln}) in this
particular case.
\end{remark}

Unlike the method used in \cite{GS94} which was based on Nevanlinna's value distribution theory, our method is
different, which is based on special functions, their asymptotic expansions, and the analytic continuation formulae for $S_{\mu,\,\nu}(\zeta)$ (\S 3.2 and \S 3.3). A crucial step in our proof is to apply the Lommel transformation (\S 2) to transform the equation (\ref{E:Chiang-Yu}) into the equation

\begin{equation}
\label{E:generalized-Lommel} \zeta^2y''(\zeta)+\zeta
y'(\zeta)+(\zeta^2-\nu^2)y(\zeta)=\sum_{j=1}^n\sigma_j
\zeta^{\mu_j+1}.
\end{equation}\smallskip

\noindent We recall from the basic differential equations theory
that the general solution of (\ref{E:generalized-Lommel}) is the
sum of {\it complementary functions}, which are the Bessel
functions, and a {\it particular integral} where each of these
particular integrals satisfies (\ref{E:generalized-Lommel}).
Since a particular integral is generally not uniquely determined, so the novelty here
is to apply standard asymptotic expansion theory (see \cite{Ol97})
to each branch of $S_{\mu,\, \nu}(\zeta)$ in order to identify
the ones that we need are precisely the particular integral whose
modulus decrease to zero and degenerate when $|\zeta| \to +\infty$, except
perhaps on the negative real axis of the $\zeta$-plane. It turns
out that these are precisely the classical Lommel function
$S_{\mu,\,\nu}(\zeta)$. We then use the inverse Lommel
transformation to (\ref{E:generalized-Lommel}) and to recover the
results for (\ref{E:Chiang-Yu}).\smallskip

The crux of the matter lies in the proof of Theorem \ref{T:Chiang-Yu}, which establishes the fact that the
function $S_{\mu,\,\nu}(L {\rm e}^{M z})$ is subnormal if and only if either $\mu+\nu$ or $\mu-\nu$ is equal to an odd positive integer.

\begin{theorem}\label{T:lommel} {Let $S_{\mu,\,\nu}(\zeta)$ be a Lommel function of an arbitrary branch. Then the entire function
$S_{\mu,\,\nu}({\rm e}^z)$ is of finite order of growth if and only
if either $\mu + \nu$ or $\mu-\nu$ is an odd positive integer $2p+1$  for some non-negative integer $p$.
In particular, the entire function $S_{\mu,\,\nu}({\rm e}^z)$
degenerates into the form {\rm (\ref{E:special-lommel-sol})} with
$\zeta={\rm e}^z$}.
\end{theorem}

 Thus,
Theorem \ref{T:Chiang-Yu} characterizes the subnormal solutions
found by Gundersen and Steinbart in \cite{GS94} to be those that correspond to the
vanishing of the complementary functions and the reduction of
particular integrals (Lommel's functions) to polynomials in $\zeta$
and $1/\zeta$.\smallskip

This paper is organized as follows. We introduce
the Lommel transformation in \S2, and show how to apply it to the equation
(\ref{E:Chiang-Yu}) to prove Theorem \ref{T:Chiang-Yu}(a). In \S 3, the
definitions of the Lommel functions are given. In particular, we
derive several new analytic continuation formulae for the Lommel
function $S_{\mu,\,\nu}(\zeta)$ in terms of the Bessel functions
of the third kind (\textit{i.e.}, Hankel functions) in the independent variable
and also in terms of the \textit{Chebyshev} polynomials of the second kind
$U_m(\cos \zeta)$ but in the parameters. In \S 4, by applying the asymptotic expansions of
$H_{\nu}^{(1)}(\zeta),\,H_{\nu}^{(2)}(\zeta)$ and
$S_{\mu,\,\nu}(\zeta)$, we can show that several entire functions
are not subnormal. The crux of the matter in the proof of Theorem
\ref{T:Chiang-Yu} that utilizes Theorem \ref{T:lommel} is to show that the function $S_{\mu,\,\nu}(L {\rm e}^{M z})$ is subnormal if and only if either $\mu+\nu$ or $\mu-\nu$ must be an odd positive integer. Although the $S_{\mu,\,\nu}(L {\rm e}^{M z})$ is single-valued, the $S_{\mu,\,\nu}(\zeta)$ so defined is multi-valued in general. The nature of the problem forces us to consider our problem for the Lommel functions for all branches. Unfortunately, we have found that the literature on the analytic continuation of the Lommel functions is inadequate, so that we have derived these new formulae in \S 3. The details of the proof of Theorem \ref{T:Chiang-Yu} is also given in \S 4. In \S 5, we prove
Theorem \ref{T:lommel} and a consequence of it will be given. In \S 6, we establish some
analogues of now classical `Quantization-type' theorems
(see e.g. \cite{La93}, Theorem 5.22) for non-homogeneous
equations as a corollary to Theorem \ref{T:Chiang-Yu}. These theorems
can be regarded as semi-classical quantitization-type results from quantum
mechanics view-point. They are followed by corresponding examples. It is here
that we identify the polynomials in (\ref{E:Gundersen-Steinbart-soln}) correspond to a
number of special polynomials: Neumann's polynomials, Gegenbauer's
polynomials and Schl\"afli's polynomials. The details of these
polynomials and Struve's functions are given in Appendix \ref{Adx:Special}.  The analytic continuation formulae and asymptotic expansions of Bessel functions are listed in Appendix \ref{Adx:Bessel}.
\medskip

\section{\bf The Lommel transformations and a proof of Theorem {\rm \ref{T:Chiang-Yu}(a)}}\label{S:lomme-transformations}
Lommel investigated transformations that involve Bessel equations
\cite{Lo1868} in 1868. Our standard references are
\cite{Wa44}, \S4.31 and \cite{Bateman}, p. 13. We mentioned that the
same transformations were also considered independently by K. Pearson (see \cite{Wa44}, p. 98) in 1880. Lommel considered the transformation

\begin{equation}
\label{E:Lommel-transformation} \zeta=\alpha x^\beta, \quad y(\zeta)=x^{\gamma}u(x),
\end{equation}\smallskip

\noindent where $x$ and $u(x)$ are the new independent and dependent
variables respectively, $\alpha, \beta \in \mathbf{C}\setminus\{0\}$
and $\gamma \in \mathbf{C}$.\smallskip

We apply this transformation to equation (\ref{E:generalized-Lommel}). It is straightforward to verify that
the function $u$ satisfies the equation

\begin{align}
\label{E:Reduced} x^2u''(x)+(2\gamma+1)x
u'(x)&+\big[\alpha^2\beta^2x^{2\beta}+(\gamma^2-\nu^2\beta^2)\big]u(x)
\notag\\
&=\sum_{j=1}^n \sigma_j \alpha^{\mu_j+1}\beta^2
x^{\beta(\mu_j+1)-\gamma}
\end{align}\smallskip

\noindent which has $x^{-\gamma} y(\alpha x^\beta)$ as its general
solution. Following the idea in \cite{CI02(2)}, we now apply a
further change of variable by
$$
x={\rm e}^{mz}, \quad f(z)=u(x)
$$ to (\ref{E:Reduced}), then we have

\begin{align}\label{E:general}
f''+2\gamma mf'+m^2\big[\alpha^2\beta^2{\rm e}^{2m\beta
z} &+(\gamma^2-\nu^2\beta^2)\big]f \notag\\
&=m^2\sum_{j=1}^n \sigma_j \alpha^{\mu_j+1}\beta^2 {\rm
e}^{m[\beta(\mu_j+1)-\gamma]z}.
\end{align}\smallskip

\noindent Choosing $m=1$ in (\ref{E:general}) and then replacing $\alpha,\beta$ and $\gamma$ by $L,M$ and $N$ respectively, yields (\ref{E:Chiang-Yu}). As we have noted in \S 1 that the general solution of (\ref{E:lommel}) is given by  a combination of the Bessel functions of first and second kinds and the Lommel function $S_{\mu,\,\nu}(\zeta)$ (see \cite{Bateman}, \S7.5.5), hence the general solution to (\ref{E:generalized-Lommel}) is

\begin{equation}
\label{E:general-soln-bessel} y(\zeta)=A J_\nu (\zeta)+B Y_\nu
(\zeta)+\sum_{j=1}^n \sigma_j S_{\mu_j,\, \nu}(\zeta).
\end{equation}\smallskip

\noindent Then the general solution of the differential equation
(\ref{E:Chiang-Yu}) is therefore given by
(\ref{E:Chiang-Yu-soln-1}). This proves Theorem
\ref{T:Chiang-Yu}(a).
\medskip

\section{\bf New formulae to the Lommel functions}
\label{S:lommel} The major part of the proof of our main theorems
consists of studying the growth of the composite function
$S_{\mu,\, \nu}({\rm e}^z)$ in the complex $z$-plane. However, the
growth of $S_{\mu,\, \nu}({\rm e}^z)$ as a subnormal solution must
be \textit{independent of the different branches} of $S_{\mu,\,
\nu}(\zeta)$, so one needs to consider its growth in all such
branches. We first note that the Lommel functions have a rather complicated
definition with respect to different subscripts (in four
different cases) even in the principal branch of $\zeta$. Since we
cannot find such analytic continuation formulae for the Lommel
functions in the literature in general, and the formulae with
respect to the different \textit{singular} subscripts in
particular (\textit{i.e.,} $\mu\pm\nu$ equals an odd negative integer), so
 we shall derive these new continuation formulae in this
section. Due to the complicated nature of the Lommel functions with respect
to different subscripts in the principal branch, so we make no apology to list some
known properties in \S \ref{Sub:Lommel-1} below before we derive
new analytic continuation formulae in \S \ref{Sub:Lommel-2} and
\S\ref{Sub:Lommel-3}. Here our main reference for the Lommel
functions are Watson \cite{Wa44}, \S10.7-10.75 {and Lommel \cite{Lo1876}}.
\medskip

\subsection{The definitions of Lommel's functions $s_{\mu,\,\nu}(\zeta)$ and $S_{\mu,\,\nu}(\zeta)$}
\label{Sub:Lommel-1} Suppose that $\mu$ and $\nu$ are complex
numbers such that none of $\mu + \nu$ and $\mu - \nu$ is an odd
negative integer. Standard variation of parameters method applied
to the equation (\ref{E:lommel}) yields a particular solution

\begin{equation*}
    \label{E:variation-formula}
    \D s_{\mu,\,\nu}(\zeta)=\frac{\pi}{2}\left[Y_\nu(\zeta) \int^\zeta t^\mu
    J_\nu(t)\, dt-J_\nu(\zeta) \int^\zeta t^\mu Y_\nu(t)\,dt \right],
\end{equation*}\smallskip

\noindent which gives (see \cite{Wa44}, \S10.7) raise to the
following expansion

\begin{align}
\label{E:small-lommel}
    s_{\mu,\,\nu}(\zeta)&=\sum_{m=0}^{+\infty}
    \frac{(-1)^m\zeta^{\mu+1+2m}}{[(\mu+1)^2-\nu^2][(\mu+3)^2-\nu^2]
    \cdots[(\mu+2m+1)^2-\nu^2]}\notag\\
    &=\frac{\zeta^{\mu+1}}{(\mu-\nu+1)(\mu+\nu+1)}\,\cdot {_1F_2}
    \left(1;\frac{1}{2}\mu-\frac{1}{2}\nu+\frac{3}{2},\frac{1}{2}\mu
    +\frac{1}{2}\nu+\frac{3}{2};-\frac{1}{4}\zeta^2\right).
\end{align}\smallskip

\noindent The above series, which begins with the term $\zeta^{\mu+1}$, is convergent for
all $\zeta$, provided that none of the parameters $\mu+\nu$ and $\mu-\nu$ is allowed to be an
odd negative integer (see \cite{Wa44}, \S10.7). Otherwise the second and
third arguments in the ${_1F_2}$ in (\ref{E:small-lommel}) would
be meaningless. This explains the restriction on $\mu\pm
\nu$ given above (see \cite{Wa44}, \S10.7). Following Watson, we
define another particular solution $S_{\mu,\,\nu}(\zeta)$  for the
equation (\ref{E:lommel}), which is also called a \textit{Lommel
function}. It is related to the $s_{\mu,\,\nu}(\zeta)$ by the
following formulae

\begin{align}
    S_{\mu,\,\nu}(\zeta)&:=s_{\mu,\,\nu}(\zeta)+ \frac{K}{\sin\nu\pi}
\bigg[J_{-\nu}(\zeta)\cos\big(\frac{\mu-\nu}{2}\pi\big)
    -J_\nu(\zeta)\cos\big(\frac{\mu+\nu}{2}\pi\big)\bigg],\label{E:lommel-lommel-1}\\
   &:=s_{\mu,\,\nu}(\zeta)+ K
\bigg[J_{\nu}(\zeta)\sin\big(\frac{\mu-\nu}{2}\pi\big)
    -Y_\nu(\zeta)\cos\big(\frac{\mu-\nu}{2}\pi\big)\bigg],\label{E:lommel-lommel-2}
\end{align}\smallskip

\noindent where

\begin{equation}
\label{E:K-constant}
K=2^{\mu-1}\Gamma\left(\frac{1}{2}\mu-\frac{1}{2}\nu+\frac{1}{2}\right)\Gamma
    \left(\frac{1}{2}\mu+\frac{1}{2}\nu+\frac{1}{2}\right),
\end{equation} \smallskip

\noindent $\mu,\, \nu \in \mathbf{C}$. The first definition
(\ref{E:lommel-lommel-1}) holds in all cases of $\mu,\, \nu$ except
 when $\nu$ is an integer. The equivalent form
(\ref{E:lommel-lommel-2}) holds even when $\nu$ is an integer. So
we shall adopt the second form (\ref{E:lommel-lommel-2}) as the
general definition for the function $S_{\mu,\,\nu}(\zeta)$
\cite{Wa44}, p. 347.\smallskip

This second Lommel function $S_{\mu,\,\nu}(\zeta)$ so defined has the advantage that it is still meaningful even when either $\mu+\nu$ or $\mu-\nu$ is an odd negative integer (see below), while $s_{\mu,\,\nu}(\zeta)$ remains undefined in (\ref{E:lommel-lommel-2}) for either of these parameter values. Since our solution to the main Theorems will involve the Lommel functions valid for \textit{all complex subscripts}, so we shall seek a way to define $S_{\mu,\,\nu}(\zeta)$ when either $\mu+\nu$ or $\mu-\nu$ is an odd negative integer. Since one easily see from the formulae (\ref{E:small-lommel}), (\ref{E:lommel-lommel-1}) and (\ref{E:K-constant}) that $S_{\mu,\,\nu}(\zeta)$ is an even function of $\nu$, it would be sufficient to derive a formula for $S_{\mu,\,\nu}(\zeta)$ when $\mu-\nu$ is an odd
negative integer and the way we define it is shown below.\smallskip

 It is known that both the $s_{\mu,\,\nu}(\zeta)$ and $S_{\mu,\,\nu}(\zeta)$ satisfy the same recurrence relation
\cite{Wa44}, \S10.72 (1) and (6):

\begin{equation}
    \label{E:lommel-at-odd}
    S_{\mu+2,\,\nu}(\zeta)=\zeta^{\mu+1}-[(\mu+1)^2-\nu^2]S_{\mu,\,\nu}(\zeta).
\end{equation}\smallskip

\noindent In particular, when $\mu\pm\nu$ is an odd negative integer,
then one may use this recurrence relation repeatedly to define
$S_{\mu,\,\nu}(\zeta)$. Indeed letting $\mu=\nu-2p-1$ and applying (\ref{E:lommel-at-odd}) as in
\cite{Wa44}, \S10.73 (1) repeatedly yields

\begin{equation}\label{E:lommel-negative}
S_{\nu-2p-1,\,
\nu}(\zeta):=\sum_{m=0}^{p-1} {(-1)^m\zeta^{\nu-2p+2m}\over
2^{2m+2}(-p)_{m+1}(\nu-p)_{m+1}} +{(-1)^pS_{\nu-1,\,
\nu}(\zeta)\over 2^{2p}p!(1-\nu)_p},
\end{equation}\smallskip

\noindent where $p$ is a positive integer. Thus the formula (\ref{E:lommel-negative}) indicates that it remains
to define the $S_{\nu-1,\, \nu}(\zeta)$. We distinguish three cases as follows:

\begin{enumerate}
  \item[(i)] If $-\nu \not \in \{0,\,1,\,2,\ldots\}$, then we apply L' Hospital's theorem \textit{once} to the relation (\ref{E:lommel-at-odd}) and obtain (\cite{Wa44}, \S10.73 (3)) that

      \begin{align}\label{E:lommel-(nu-1)-original}
      S_{\nu-1,\,\nu}(\zeta)&={1\over 4}\zeta^\nu\Gamma(\nu)\sum_{m=0}^{+\infty} \frac{(-1)^m({1\over 2}\zeta)^{2m}}{m! \Gamma(\nu+m+1)} \big[2\log\zeta-A(m)\big] \\
      &\quad-2^{\nu-2}\pi\Gamma(\nu) Y_\nu(\zeta),\notag
      \end{align}\smallskip

\noindent where $A(m)=2\log 2+\psi(\nu+m+1)+\psi(m+1)$ and $\psi$ is the digamma function. For our later applications, let us rewrite the formula (\ref{E:lommel-(nu-1)-original}) as

\begin{align}\label{E:lommel-(nu-1)}
S_{\nu-1,\,\nu}(\zeta)&=\Gamma(\nu)\Bigg[2^{\nu-1} J_\nu(\zeta)\log\zeta-2^{\nu-2}\pi
Y_\nu(\zeta)\\
&\quad\quad\quad -\frac{\zeta^\nu}{4}\sum_{m=0}^{+\infty}
\frac{(-1)^m(\frac{1}{2}\zeta)^{2m}A(m)}{m!\Gamma(\nu+m+1)}\Bigg].\notag
\end{align} \smallskip

  \item[(ii)] If $\nu=0$, then we apply L' Hospital's Theorem \textit{twice} to the formula

$$S_{\mu,\,0}(\zeta)=\frac{\zeta^{\mu+1}-S_{\mu+2,\,0}(\zeta)}{(\mu+1)^2}$$\smallskip

\noindent to obtain (\cite{Wa44}, \S10.73 (4))

\begin{align}\label{E:lommel-negative-one-definition}
S_{-1,\,0}(\zeta)&=\frac{1}{2}\sum_{m=0}^{+\infty}
\frac{(-1)^m\big(\frac{1}{2}\zeta\big)^{2m}}{(m!)^2}\\
&\quad\times \bigg\{\bigg[\log
\frac{\zeta}{2}-\psi(m+1)\bigg]^2-\frac{1}{2}\psi'(m+1)+\frac{\pi^2}{4}\bigg\}. \notag
\end{align} \smallskip

\noindent Again for the easy of our applications later, let us rewrite the formula (\ref{E:lommel-negative-one-definition}) as

\begin{align}\label{E:lommel-negative-one}
S_{-1,\,0}(\zeta)&=\frac{1}{2}J_0(\zeta)(\log
\zeta)^2+\frac{1}{2}\sum_{m=0}^{+\infty}\frac{(-1)^m(\frac{1}{2}\zeta)^{2m}B(m)}{(m!)^2}\\
&\quad-\Bigg[J_0(\zeta)\log 2+\sum_{m=0}^{+\infty}\frac{(-1)^m(\frac{1}{2}\zeta)^{2m}}{(m!)^2}\psi(m+1)\Bigg]\log \zeta, \notag
\end{align}\smallskip

\noindent where $\D B(m)=\bigg[\log 2+\psi(m+1)\bigg]^2-\frac{1}{2}\psi'(m+1)+\frac{\pi^2}{4}$.

  \item[(iii)] If $\nu=-n$ for a positive integer $n$, then we can apply the formula (see \cite{Wa44}, \S10.72 (7))

  \begin{equation}\label{E:lommel-derivative}
  S_{\mu,\,\nu}'(\zeta) +\frac{\nu}{\zeta}S_{\mu,\,\nu}(\zeta)=(\mu+\nu-1)S_{\mu-1,\,\nu-1}(\zeta)
  \end{equation}\smallskip

  \noindent to obtain the formula (\cite{Lo1876}, Eqn. (XIX), p. 439)

  \begin{equation}\label{E:lommel-derivative-n}
  S_{-n-1,\,-n}(\zeta)=\frac{(-1)^n\zeta^n}{n!}\cdot\frac{{\rm d}^n}{{\rm d}(\zeta^2)^n}(S_{-1,\,0}(\zeta)).
  \end{equation}\smallskip

  \noindent Hence the function $S_{-n-2p-1,\,-n}(\zeta)$ can be defined from (\ref{E:lommel-negative}), (\ref{E:lommel-negative-one}) and (\ref{E:lommel-derivative-n}).\smallskip

  {Since the $S_{-1,\,0}(\zeta)$ is a solution of the differential equation (\ref{E:lommel}) with $\mu=-1$ and $\nu=0$, we deduce from the formula (\ref{E:lommel-derivative-n}) that}

  \begin{align}\label{E:lommel-negative-integer-general}
  S_{-n-1,\,-n}(\zeta)=\frac{(-1)^n\zeta^{-n}}{2^nn!}\big[A_n(\zeta)+B_n(\zeta)S_{-1,\,0}(\zeta)+\zeta C_n(\zeta)S_{-1,\,0}'(\zeta)\big],
  \end{align}\smallskip

  \noindent where $A_n(\zeta),\, B_n(\zeta)$ and $C_n(\zeta)$ are polynomials in $\zeta$ of degree at most $n$ such that $A_1(\zeta)=B_1(\zeta) \equiv 0,\,C_1(\zeta) \equiv 1$, and when $n \ge 2$, that they satisfy the following recurrence relations:

  \begin{equation}\label{E:recurrence-A_n-B_n-C_n}
  \begin{split}
  A_n(\zeta)&=-2(n-1)A_{n-1}(\zeta)+\zeta A_{n-1}'(\zeta)+C_{n-1}(\zeta),\\
  B_n(\zeta)&=-2(n-1)B_{n-1}(\zeta)+\zeta B_{n-1}'(\zeta)-\zeta^2 C_{n-1}(\zeta),\\
  C_n(\zeta)&=-2(n-1)C_{n-1}(\zeta)+B_{n-1}(\zeta)+\zeta C_{n-1}'(\zeta).
  \end{split}
  \end{equation}\smallskip
\end{enumerate}

\begin{remark}
{\rm We remark about a property of the formula (\ref{E:lommel-negative-integer-general}) that will be used in the proof of Proposition 4.4 as follows: Let $P(1/\zeta)=a_0+a_1/\zeta+\cdots+a_n/\zeta^n$ be a polynomial in $1/\zeta$ of degree $n$. Then we can show by the method of comparing coefficients that $P(1/\zeta)$ does not satisfy the differential equation (\ref{E:lommel}) when $\mu=-n-1$ and $\nu=-n$. This shows that we cannot have $B_n(\zeta)S_{-1,\,0}(\zeta)+\zeta C_n(\zeta)S_{-1,\,0}'(\zeta) \equiv 0$ for any positive integer $n$ in (\ref{E:lommel-negative-integer-general}).}
\end{remark}

\begin{remark}\rm Since it is {not} straightforward to derive an analytic continuation formula from
the formula (\ref{E:lommel-negative-one}), we need to further
simplify it. In fact, we have

\begin{equation*}\label{E:bessel-second-zero}
\pi Y_0(\zeta)=2\Bigg[J_0(\zeta)\log\frac{\zeta}{2}-
\sum_{m=0}^{+\infty}\frac{(-1)^m(\frac{1}{2}\zeta)^{2m}}{(m!)^2}\psi(m+1)\Bigg]
\end{equation*}\smallskip

\noindent (see \cite{Wa44}, Eqn. (2), p. 60 and Eqn. (2), p. 64), so
the equation (\ref{E:lommel-negative-one}) can be further written
as

\begin{align}\label{E:lommel-negative-one-further}
S_{-1,\,0}(\zeta) =-\frac{1}{2}J_0(\zeta)(\log
\zeta)^2+\frac{\pi}{2}Y_0(\zeta)\log
\zeta+\frac{1}{2}\sum_{m=0}^{+\infty}\frac{(-1)^m(\frac{1}{2}\zeta)^{2m}B(m)}{(m!)^2}.
\end{align}\smallskip
\end{remark}

\noindent We are ready to derive new analytic continuation formulae for the Lommel functions first with respect to regular subscripts in \S3.2 (Lemma \ref{L:lommel-continuation-hankel} and Theorem \ref{T:lommel-continuation-general}), and then with respect to the singular subscripts in \S3.3 (Lemmae \ref{L:lommel-continuation-general-1} to \ref{L:lommel-continuation-odd-negative-general}).

\subsection{New analytic continuation formula for $S_{\mu,\,\nu}(\zeta)$ when none of $\mu \pm \nu$ is an odd negative integer} \label{Sub:Lommel-2}

The proof of the Theorem \ref{T:Chiang-Yu}(b) relies heavily on the fact
that the function $S_{\mu,\,\nu}({\rm e}^z)$ so defined is
independent of a particular branch of $S_{\mu,\,\nu}(\zeta)$ under
consideration. We thus need to derive analytic continuation formulae for
all such branches of the function. We first derive the continuation formula with a full proof in Lemma
\ref{L:lommel-continuation-hankel} below.

\begin{lemma}\label{L:lommel-continuation-hankel} We have

\begin{equation*}\label{E:lommel-continuation-hankel} S_{\mu,\,\nu}(\zeta {\rm e}^{-\pi i})
= -{\rm e}^{-\mu \pi i}S_{\mu,\,\nu}(\zeta)+K_+
H_{\nu}^{(1)}(\zeta),
\end{equation*}\smallskip

\noindent where $\D K_+=Ki[1+{\rm e}^{(-\mu+\nu)\pi
i}]\cos\bigl(\frac{\mu+\nu}{2}\pi\bigr)$, and $K$ is given  by
{\rm (\ref{E:K-constant})}.
\end{lemma}
\medskip

\begin{proof} Let $m \in {\bf Z}$, then it is easy to check that

\begin{equation}\label{E:small-lommel-analytic-continuation}
s_{\mu,\,\nu}(\zeta {\rm e}^{m\pi i})=(-1)^m{\rm e}^{m\mu\pi
i}s_{\mu,\,\nu}(\zeta)
\end{equation}\smallskip

\noindent holds. Let $K$ be given by (\ref{E:K-constant}). It follows from (\ref{E:lommel-lommel-2}), (\ref{E:small-lommel-analytic-continuation}), (\ref{E:bessel-first-continuation}) and
(\ref{E:bessel-second-continuation}) with $m=-1$ that

\begin{align*}
S_{\mu,\,\nu}(\zeta {\rm e}^{-\pi i})&= s_{\mu,\,\nu}(\zeta {\rm
e}^{-\pi i})+K\bigg[\sin\bigl({\mu-\nu\over 2}\pi\bigr)
J_{\nu}(\zeta {\rm e}^{-\pi i})- \cos\bigl({\mu-\nu\over
2}\pi\bigr)Y_\nu(\zeta {\rm e}^{-\pi i})\bigg]\cr
                         &= -{\rm e}^{-\mu\pi i}s_{\mu,\,\nu}(\zeta)+
K\times \bigg\{{\rm e}^{-\nu\pi i}\sin\bigl({\mu-\nu\over
2}\pi\bigr) J_{\nu}(\zeta)\cr
                         &\qquad- \cos\bigl({\mu-\nu\over
2}\pi\bigr) \bigl[{\rm e}^{\nu\pi i}Y_\nu(\zeta)-2i\cos(\nu\pi)
J_\nu(\zeta)\bigr]\bigg\}.
\end{align*}\smallskip

\noindent We now substitute for $s_{\mu,\,\nu}(\zeta)$ in terms of
$S_{\mu,\,\nu}(\zeta)$ from (\ref{E:lommel-lommel-2}) in the above
equation to yield an analytic continuation formula so that

\begin{align}\label{E:lommel-continuation-bessel} S_{\mu,\,\nu}(\zeta {\rm e}^{-\pi
i})&= -{\rm e}^{-\mu\pi
i}S_{\mu,\,\nu}(\zeta)+K\bigg\{\bigg[\bigl({\rm e}^{-\mu\pi
i}+{\rm e}^{-\nu\pi i}\bigr)\sin\bigl({\mu-\nu\over
2}\pi\bigr)\bigg]J_\nu(\zeta)\\
                         &\hspace{2.5cm}+2i\cos\bigl({\mu-\nu\over
                         2}\pi\bigr)\cos(\nu\pi)
J_\nu(\zeta)\bigg\} \notag\\
&\quad -K({\rm e}^{-\mu\pi i}+{\rm e}^{\nu\pi
i}\bigr)\cos\bigl({\mu-\nu\over 2}\pi\bigr)Y_\nu(\zeta).\notag
\end{align}\smallskip

\noindent Replacing the Bessel functions of the first and
second kinds in (\ref{E:lommel-continuation-bessel}) by the Hankel
functions (\ref{E:hankels-definition}) yields

\begin{equation*}
S_{\mu,\,\nu}(\zeta {\rm e}^{-\pi i})= -{\rm e}^{-\mu\pi
i}S_{\mu,\,\nu}(\zeta)+K_+H_\nu^{(1)}(\zeta)+K_-H_\nu^{(2)}(\zeta),
\end{equation*}\smallskip

\noindent where

\begin{align*}
K_{\pm}&={K\over 2}\bigg[({\rm e}^{-\mu\pi i}+{\rm e}^{-\nu\pi
i})\sin\bigl({\mu-\nu\over 2}\pi\bigr)+2i\cos\bigl({\mu-\nu\over
2}\pi\bigr)\cos(\nu\pi)\\
&\quad \pm i({\rm e}^{-\mu\pi i}+{\rm e}^{\nu\pi
i})\cos\bigl({\mu-\nu\over 2}\pi\bigr)\bigg].
\end{align*}\smallskip

Substituting $\D \cos \theta={{\rm e}^{\theta i}+{\rm e}^{-\theta
i} \over 2}$ and $\D \sin \theta={{\rm e}^{\theta i}-{\rm e}^{-\theta
i} \over 2i}$ into the above equations yields

$$K_+=Ki[1+{\rm e}^{(-\mu+\nu)\pi
i}]\cos\bigl(\frac{\mu+\nu}{2}\pi\bigr)\quad \mbox{and}\quad
K_-=0.$$\smallskip
\end{proof}
\smallskip

Let us write for each integer $m$ that

\begin{equation}\label{E:small-polynomial} U_{m-1}(\cos \zeta):=\frac{\sin m\zeta}{\sin
\zeta}
\end{equation}\smallskip

\noindent which is the \textit{Chebyshev polynomials of the second kind}, see \cite{Bateman}. It is a polynomial
of $\cos \zeta$ of degree $(m-1)$. We note the elementary facts that

$$U_0(\cos \nu\pi)=1,\ U_1(0)=0,\ U_{m-1}(1)=\lim_{\zeta \to 0}\frac{\sin m\zeta}{\sin \zeta}=m$$\smallskip

\noindent and

\begin{equation}\label{E:p_m}
    \D U_{m-1}(\cos k\pi)=\lim_{\nu \to k}\frac{\sin m\nu\pi}{\sin \nu\pi}=m(-1)^{k(m-1)}
\end{equation}\smallskip

\noindent hold for each integer $k$. Thus we have the following analytic continuation formula:

\begin{theorem}\label{T:lommel-continuation-general} Let $m$ be a non-zero integer, and $\mu\pm \nu\not= 2p+1$ for any integer $p$.

\begin{enumerate}
  \item[(a)] We have

\begin{align}\label{E:lommel-continuation-general}
 S_{\mu,\,\nu}(\zeta {\rm e}^{-m\pi i})&=(-1)^m {\rm e}^{-m\mu\pi i}
 S_{\mu,\,\nu}(\zeta)+K_+\big[P_m(\cos\nu\pi,{\rm e}^{-\mu \pi i})H_\nu^{(1)}(\zeta)\\
 &\quad+{\rm e}^{-\nu\pi i }Q_m(\cos\nu\pi,{\rm e}^{-\mu \pi i})H_\nu^{(2)}(\zeta)\big],\notag
\end{align}\smallskip

\noindent where $K_+$ is given in the Lemma {\rm
\ref{L:lommel-continuation-hankel}}, $P_m(\cos\nu\pi,{\rm e}^{-\mu
\pi i})$ and $Q_m(\cos\nu\pi,{\rm e}^{-\mu \pi i})$ are rational functions
of $\cos \nu\pi$ and ${\rm e}^{-\mu\pi i}$ given by

\begin{align}
P_m&(\cos \nu\pi,\,{\rm e}^{-\mu\pi i}) \notag\\
&=\frac{U_{m-1}(\cos \nu\pi) +{\rm e}^{-\mu\pi i}U_m(\cos \nu\pi)+(-1)^{m+1}{\rm e}^{-(m+1)\mu\pi i}}{[1+{\rm e}^{-(\mu+\nu)\pi i}][1+{\rm e}^{-(\mu-\nu)\pi i}]} \label{E:coefficient-P}
\end{align}\smallskip

\noindent and

\begin{align}
Q_m&(\cos \nu\pi,\,{\rm e}^{-\mu\pi i}) \notag\\
&=\frac{U_{m-2}(\cos \nu\pi) +{\rm e}^{-\mu\pi i}U_{m-1}(\cos \nu\pi)+(-1)^m{\rm e}^{-m\mu\pi i}}{[1+{\rm e}^{-(\mu+\nu)\pi i}][1+{\rm e}^{-(\mu-\nu)\pi i}]} \label{E:coefficient-Q},
\end{align}\smallskip

\noindent where $U_m(\cos \nu\pi)$ is given by {\rm (\ref{E:small-polynomial})}.\smallskip

  \item[(b)]  Furthermore, the coefficients $P_m(\cos\nu\pi,\, {\rm e}^{-\mu \pi i})$ and $Q_m(\cos\nu\pi,\, {\rm e}^{-\mu \pi i})$ are not identically zero \textit{simultaneously} for all $\mu,\nu$ and all non-zero integers $m$.
\end{enumerate}
\end{theorem}\medskip

Before proving Theorem \ref{T:lommel-continuation-general}, we need the following {relations} which can be derived easily from the definition (\ref{E:small-polynomial}): For any integer $m$, we have

\begin{equation}\label{E:identities}
\begin{split}
&U_{-m}(\cos \nu\pi)=-U_{m-2}(\cos \nu\pi),\\
&U_{-m-2}(\cos \nu\pi)=-U_m(\cos \nu\pi),\\
&U_{m-1}^2(\cos \nu\pi)+U_m(\cos \nu\pi)U_{-m}(\cos \nu\pi)=1.
\end{split}
\end{equation}\smallskip

\begin{proof}[Proof of Theorem \ref{T:lommel-continuation-general}] In fact, we first prove the claim that the formula (\ref{E:lommel-continuation-general})
holds with the expressions $P_m(\cos\nu\pi,\, {\rm e}^{-\mu \pi i})$ and $Q_m(\cos\nu\pi,\, {\rm e}^{-\mu \pi i})$ given by

\begin{align}\label{E:coefficient-P-claim}
P_m&(\cos\nu\pi,\, {\rm e}^{-\mu \pi i}) \notag\\
&=\left\{
  \begin{array}{ll}
    \D \sum_{j=0}^m (-1)^j{\rm e}^{-j\mu\pi i}U_{m-j-1}(\cos\nu\pi), & \hbox{if $m > 0$}; \\
    & \\
    (-1)^{m+1}{\rm e}^{-m \mu\pi i}\big[P_{-m}(\cos\nu\pi,{\rm e}^{-\mu \pi i})U_m(\cos\nu\pi) & \\
-Q_{-m}(\cos\nu\pi,{\rm e}^{-\mu \pi i})U_{m-1}(\cos\nu\pi)\big],& \hbox{if
$m<0$,}
  \end{array}
\right.
\end{align}\smallskip

\noindent and

\begin{align}\label{E:coefficient-Q-claim}
Q_m&(\cos\nu\pi,\, {\rm e}^{-\mu \pi i})\notag\\
&=\left\{
  \begin{array}{ll}
    \D \sum_{j=0}^{m-1} (-1)^j{\rm e}^{-j\mu\pi i}U_{m-j-2}(\cos \nu\pi),& \hbox{if $m > 0$}; \\
    & \\
    (-1)^{m+1} {\rm e}^{-m\mu \pi i}
  \big[P_{-m}(\cos\nu\pi,{\rm e}^{-\mu \pi i})U_{m-1}(\cos\nu\pi) & \\
   +Q_{-m}(\cos\nu\pi,{\rm e}^{-\mu \pi i})U_{-m}(\cos\nu\pi)\big], & \hbox{if $m<0$.}
  \end{array}
\right.
\end{align}\smallskip

We apply induction on positive integers $m$. Suppose $m=1$, then Lemma \ref{L:lommel-continuation-hankel} shows
that (\ref{E:lommel-continuation-general}) holds with

$$P_1(\cos\nu\pi,{\rm e}^{-\mu \pi i})\equiv 1\quad \mbox{and}\quad
Q_1(\cos\nu\pi,{\rm e}^{-\mu \pi i})\equiv 0$$\smallskip

\noindent as given by (\ref{E:coefficient-P-claim}) and (\ref{E:coefficient-Q-claim}) respectively, which are
trivial rational functions in $\cos \nu\pi$ and ${\rm e}^{-\mu\pi i}$.
\smallskip

We note that the analytic continuation formulae
(\ref{E:bessel-third-1-continuation}) and
(\ref{E:bessel-third-2-continuation}) for $H_\nu^{(1)}(\zeta)$ and
$H_\nu^{(2)}(\zeta)$ can be rewritten as

\begin{align}
H_\nu^{(1)}(\zeta {\rm e}^{m\pi i})&=U_{-m}(\cos\nu\pi)H_\nu^{(1)}(\zeta)-{\rm e}^{-\nu\pi i}U_{m-1}(\cos\nu\pi)
H_\nu^{(2)}(\zeta),\label{E:bessel-third-1-continuation-rewritten}\\
H_\nu^{(2)}(\zeta {\rm e}^{m\pi i})&=U_{m}(\cos\nu\pi)H_\nu^{(2)}(\zeta)+{\rm e}^{\nu\pi i}U_{m-1}(\cos\nu\pi)
H_\nu^{(1)}(\zeta).
\label{E:bessel-third-2-continuation-rewritten}
\end{align}\smallskip

Let us suppose that the formula
(\ref{E:lommel-continuation-general}) holds for $m=k$, where $k\in
\mathbf{N}$, \textit{i.e.},

\begin{align*}
S_{\mu,\,\nu}(\zeta {\rm e}^{-k\pi i})&=(-1)^k {\rm e}^{-k\mu\pi
i} S_{\mu,\,\nu}(\zeta)+K_+\big[P_k(\cos\nu\pi,{\rm e}^{-\mu \pi
i})H_\nu^{(1)}(\zeta)\\
&\qquad+{\rm e}^{-\nu \pi i }Q_k(\cos\nu\pi,{\rm e}^{-\mu \pi
i})H_\nu^{(2)}(\zeta)\big].
\end{align*} \smallskip

\noindent We observe that the polynomial (\ref{E:small-polynomial})
satisfies the relation

\begin{equation}\label{E:Tchebichef-relation}
U_{-m-1}(\cos\zeta)=-U_{m-1}(\cos\zeta)
\end{equation}\smallskip

\noindent for any integer $m$. Thus we have by Lemma \ref{L:lommel-continuation-hankel}, (\ref{E:bessel-third-1-continuation-rewritten}) and the relation (\ref{E:Tchebichef-relation})\footnote{The relation (\ref{E:Tchebichef-relation}) is obtained by applying an {inductive argument} to $Q_{k+1}(\cos \nu\pi,\,{\rm e}^{-\mu\pi i})$ via $P_k(\cos \nu\pi,\,{\rm e}^{-\mu\pi i})$ and $Q_k(\cos \nu\pi,\,{\rm e}^{-\mu\pi i})$.} that,

\begin{align*}
S_{\mu,\,\nu}(\zeta {\rm e}^{-(k+1)\pi i})&=-{\rm e}^{-\mu \pi
i}S_{\mu,\,\nu}(\zeta {\rm e}^{-k\pi i})+K_+H_\nu^{(1)}(\zeta {\rm e}^{-k\pi i})\\
&=-{\rm e}^{-\mu\pi i} \big\{(-1)^k {\rm e}^{-k\mu\pi
i}S_{\mu,\,\nu}(\zeta)+K_+\big[P_k(\cos\nu\pi,{\rm e}^{-\mu \pi
i})H_\nu^{(1)}(\zeta)\\
&\hspace{2.5cm}+{\rm e}^{-\nu \pi i}Q_k(\cos\nu\pi,{\rm e}^{-\mu \pi i})H_\nu^{(2)}(\zeta)\big]\big\}\\
&\qquad+K_+\big[U_{k}(\cos\nu\pi)
H_\nu^{(1)}(\zeta)-{\rm e}^{-\nu\pi i}U_{-k-1}(\cos\nu\pi) H_\nu^{(2)}(\zeta)\big]\\
&=(-1)^{k+1} {\rm e}^{-(k+1)\mu\pi i}
S_{\mu,\,\nu}(\zeta)+K_+\big[P_{k+1}(\cos\nu\pi,{\rm e}^{-\mu \pi
i})H_\nu^{(1)}(\zeta)\\
&\qquad+{\rm e}^{-\nu\pi i }Q_{k+1}(\cos\nu\pi,{\rm e}^{-\mu \pi
i})H_\nu^{(2)}(\zeta)\big],
\end{align*}\smallskip

\noindent where $P_{k+1}(\cos\nu\pi,{\rm e}^{-\mu \pi i})$ and
$Q_{k+1}(\cos\nu\pi,{\rm e}^{-\mu \pi i})$ are expressions
matching exactly the formulae (\ref{E:coefficient-P-claim}) and
(\ref{E:coefficient-Q-claim}) ($m>0$) respectively. We conclude, by induction, that the formula (\ref{E:lommel-continuation-general}) holds for all positive integers $m$.\smallskip

For a negative integer $m$, $-m$ must be positive and then we
apply the formula (\ref{E:lommel-continuation-general}) for
positive $-m$,

\begin{align}\label{E:lommel-continuation-general2}
S_{\mu,\,\nu}(\zeta)&=(-1)^{m} {\rm e}^{-m\mu\pi i}
\big\{S_{\mu,\,\nu}(\zeta {\rm e}^{m\pi i})-K_+\big[P_{-m}(\cos\nu\pi,{\rm e}^{-\mu \pi i})H_\nu^{(1)}(\zeta)\\
&\quad+{\rm e}^{-\nu\pi i}Q_{-m}(\cos\nu\pi,{\rm e}^{-\mu \pi
i})H_\nu^{(2)}(\zeta)\big]\big\}.\notag
\end{align}\smallskip

\noindent Then we replace $\zeta$ by $\zeta {\rm e}^{-m\pi i}$ in the formula (\ref{E:lommel-continuation-general2}) to get

\begin{align}\label{E:lommel-continuation-general3}
S_{\mu,\,\nu}(\zeta {\rm e}^{-m\pi i})&=(-1)^{m} {\rm e}^{-m\mu\pi i}\big\{S_{\mu,\,\nu}(\zeta)-K_+\big[P_{-m}(\cos\nu\pi,{\rm e}^{-\mu \pi i})\times\\
&\quad H_\nu^{(1)}(\zeta {\rm e}^{-m\pi i})+{\rm e}^{-\nu\pi i}Q_{-m}(\cos\nu\pi,{\rm e}^{-\mu \pi
i})H_\nu^{(2)}(\zeta {\rm e}^{-m\pi i})\big]\big\}.\notag
\end{align}\smallskip

\noindent Thus the desired results for (\ref{E:coefficient-P-claim}) and
(\ref{E:coefficient-Q-claim}) ($m<0$) follows from the formula (\ref{E:lommel-continuation-general3}) and the continuation formulae (\ref{E:bessel-third-1-continuation-rewritten})  and
(\ref{E:bessel-third-2-continuation-rewritten}) with $\zeta$ replaced by $\zeta {\rm e}^{-m\pi i}$. This proves our claim. \smallskip

Now we show that the two formulae (\ref{E:coefficient-P-claim}) and (\ref{E:coefficient-Q-claim}) can be reduced to the formulae (\ref{E:coefficient-P}) and (\ref{E:coefficient-Q}) respectively. We note that the expression of $P_m(\cos\nu\pi,\,{\rm e}^{-\mu\pi i})$ when $m > 0$ can be simplified as follows:

\begin{align}\label{E:P_m_simplify}
P_m&(\cos\nu\pi,\,{\rm e}^{-\mu\pi i}) \notag\\
&=\sum_{j=0}^m (-1)^j{\rm e}^{-j\mu\pi i}U_{m-j-1}(\cos \nu\pi) \notag\\
&=\frac{1}{2i\sin \nu\pi}\sum_{j=0}^m (-1)^j {\rm e}^{-j\mu \pi i}\big[{\rm e}^{(m-j)\nu\pi i}-{\rm e}^{-(m-j)\nu \pi i}\big] \notag\\
&=\frac{1}{2i\sin \nu\pi} \Bigg[{\rm e}^{m\nu\pi i} \sum_{j=0}^m (-1)^j {\rm e}^{-(\mu+\nu)j\pi i}-
{\rm e}^{-m\nu\pi i} \sum_{j=0}^m (-1)^j {\rm e}^{-(\mu-\nu)j\pi i}\Bigg] \notag\\
&=\frac{{\rm e}^{m\nu\pi i}}{2i\sin\nu\pi}\sum_{j=0}^m (-1)^j {\rm e}^{-(\mu+\nu)j\pi i}-\frac{{\rm e}^{-m\nu\pi i}}{2i\sin\nu\pi}\sum_{j=0}^m (-1)^j {\rm e}^{-(\mu-\nu)j\pi i} \notag\\
&=\frac{{\rm e}^{m\nu\pi i}}{2i\sin\nu\pi} \cdot \frac{1-(-1)^{m+1}{\rm e}^{-(m+1)(\mu+\nu)\pi i}}{1+{\rm e}^{-(\mu+\nu)\pi i}} \notag\\
&\hspace{2.5cm}-\frac{{\rm e}^{-m\nu\pi i}}{2i\sin\nu\pi}\cdot \frac{1-(-1)^{m+1}{\rm e}^{-(m+1)(\mu-\nu)\pi i}}{1+{\rm e}^{-(\mu-\nu)\pi i}} \notag\\
&=\frac{1}{2i\sin\nu\pi {\rm e}^{m\nu\pi i}[1+{\rm e}^{-(\mu+\nu)\pi i}][1+{\rm e}^{-(\mu-\nu)\pi i}]}\\
&\quad \times \big\{{\rm e}^{2m\nu\pi i}[1+{\rm e}^{-(\mu-\nu)\pi i}][1+(-1)^m{\rm e}^{-(m+1)(\mu+\nu)\pi i}] \notag\\
&\hspace{2.5cm} -[1+{\rm e}^{-(\mu+\nu)\pi i}][1+(-1)^m{\rm e}^{-(m+1)(\mu-\nu)\pi i}]\big\}. \notag
\end{align}\smallskip

\noindent When we expand the products in the numerator in the equation (\ref{E:P_m_simplify}), the term $(-1)^m {\rm e}^{m\nu\pi i}{\rm e}^{-(m+2)\mu\pi i}$ is eliminated and the remaining terms can be simplified to

\begin{align*}
&\quad {\rm e}^{m\nu\pi i}({\rm e}^{m\nu\pi i}-{\rm e}^{-m\nu\pi i})+{\rm e}^{-\mu\pi i}{\rm e}^{m\nu\pi i}[{\rm e}^{(m+1)\nu\pi i}-{\rm e}^{-(m+1)\nu\pi i}]\\
&\hspace{2.5cm}+(-1)^m {\rm e}^{-(m+1)\mu\pi i}{\rm e}^{m\nu\pi i} ({\rm e}^{-\nu\pi i}-{\rm e}^{\nu\pi i})\\
&={\rm e}^{m\nu\pi i} \cdot 2i\sin m\nu\pi+{\rm e}^{-\mu\pi i}{\rm e}^{m\nu\pi i}\cdot 2i\sin(m+1)\nu\pi\\
&\hspace{2.5cm}-(-1)^m{\rm e}^{-(m+1)\mu\pi i}{\rm e}^{m\nu\pi i}\cdot 2i\sin \nu\pi\\
&=2i{\rm e}^{m\nu\pi i}\big[\sin m\nu\pi +{\rm e}^{-\mu\pi i}\sin(m+1)\nu\pi-(-1)^m{\rm e}^{-(m+1)\mu\pi i}\sin \nu\pi\big],
\end{align*}\smallskip

\noindent and thus the equation (\ref{E:P_m_simplify}) yields the formula (\ref{E:coefficient-P}) for all positive integers $m$. Since the formulae (\ref{E:coefficient-P-claim}) and (\ref{E:coefficient-Q-claim}) are connected by

\begin{align}\label{E:Q_m-equivalent-form}
Q_m(\cos \nu\pi,\,{\rm e}^{-\mu\pi i})&=-{\rm e}^{\mu\pi i}\big[P_m(\cos \nu\pi,\,{\rm e}^{-\mu\pi i})-U_{m-1}(\cos\nu \pi)\big],
\end{align}\smallskip

\noindent so when $m > 0$, the expression (\ref{E:Q_m-equivalent-form}) can be written, after applying the (\ref{E:coefficient-P}), as

\begin{align*}
Q_m(\cos \nu\pi,\,{\rm e}^{-\mu\pi i})&=\frac{1}{[1+{\rm e}^{-(\mu+\nu)\pi i}][1+{\rm e}^{-(\mu-\nu)\pi i}]}\\
&\quad \times \big[({\rm e}^{\nu\pi i}+{\rm e}^{-\nu\pi i}) U_{m-1}(\cos \nu\pi)-U_m(\cos \nu\pi)\\
&\qquad+{\rm e}^{-\mu\pi i}U_{m-1}(\cos \nu\pi)+(-1)^m {\rm e}^{-m\mu\pi i}\big].
\end{align*}\smallskip

\noindent Now it follows from the definition (\ref{E:small-polynomial}) and the compound angle formulae for sine function that the first two terms in the numerator in the above equation is exactly $U_{m-2}(\cos \nu\pi)$, thus proving that the formula (\ref{E:coefficient-Q}) holds when $m > 0$.\smallskip

However, when $m$ is a negative integer, $-m$ is a positive integer. We substitute $-m$ into the formulae (\ref{E:coefficient-P}) and (\ref{E:coefficient-Q}) respectively. Substituting the resulting $P_{-m}(\cos\nu\pi,\, {\rm e}^{-\mu\pi i})$ and $Q_{-m}(\cos\nu\pi,\, {\rm e}^{-\mu\pi i})$ into the expression (\ref{E:coefficient-P-claim}), we obtain

\begin{align*}
P_m&(\cos\nu\pi,\, {\rm e}^{-\mu\pi i})\\
&=\frac{(-1)^{m+1}{\rm e}^{-m \mu\pi i}}{[1+{\rm e}^{-(\mu+\nu)\pi i}][1+{\rm e}^{-(\mu-\nu)\pi i}]}\\
&\quad \times \big[U_{-m-1}(\cos\nu\pi)U_m(\cos\nu\pi)-U_{-m-2}(\cos\nu\pi)U_{m-1}(\cos\nu\pi)\\
&\qquad+{\rm e}^{-\mu \pi i}U_{-m}(\cos\nu\pi)U_m(\cos\nu\pi)-{\rm e}^{-\mu \pi i}U_{-m-1}(\cos\nu\pi)U_{m-1}(\cos\nu\pi)\\
&\qquad+(-1)^{-m+1} {\rm e}^{(m-1)\mu \pi i}U_{m}(\cos\nu\pi)+(-1)^{-m+1} {\rm e}^{m\mu \pi i}U_{m-1}(\cos\nu\pi)\big].
\end{align*}\smallskip

\noindent The identities in (\ref{E:identities}) imply that the first two terms and the following two terms in the numerator in the above equation vanish { identically} and equal to ${\rm e}^{-\mu \pi i}$, respectively. But this is exactly the formula (\ref{E:coefficient-P}) in the case $m<0$. The validity of (\ref{E:coefficient-Q}), when $m<0$, can be obtained similarly. This completes the proof of (a). \smallskip

In order to prove (b), we note that the result of (a) implies that the expression (\ref{E:Q_m-equivalent-form}) holds for all non-zero integers $m$. Thus it is easy to deduce from (\ref{E:Q_m-equivalent-form}) that it is impossible for $P_m(\cos\nu\pi,\,{\rm e}^{-\mu\pi i})$ and $Q_m(\cos\nu\pi,\,{\rm e}^{-\mu\pi i})$ to be identically zero \textit{simultaneously} for all non-zero integers $m$, thus completing the proof of the Theorem.
\end{proof}

We note that when $\mu\pm \nu=2p+1$ where $p$ is a non-negative integer, then the constant $K_+=0$ . In fact, according to the Lemma \ref{L:lommel-continuation-hankel} it is easy to see that the continuation formula is given trivially in Remark \ref{R:lommel-terminate}.
\smallskip

\subsection{New analytic continuation formulae for $S_{\mu,\,\nu}(\zeta)$ when either $\mu+\nu$ or $\mu-\nu$ is an odd negative integer}\label{Sub:Lommel-3}

 We recall that in this case that the $s_{\mu,\,\nu}(\zeta)$ is undefined, so we cannot apply the
definition (\ref{E:lommel-lommel-2}). Instead, we shall use
the formulae (\ref{E:lommel-(nu-1)}), (\ref{E:lommel-negative-integer-general}) and (\ref{E:lommel-negative-one-further}) to obtain analytic continuation formulae for $S_{\nu-1,\,\nu}(\zeta)$ and then for $S_{\nu-2p-1,\,\nu}(\zeta)$, where $p$ is a non-negative integer. We shall only deal with the case $\mu-\nu=-2p-1$ in the following argument, while the remaining case $\mu+\nu=-2p-1$ can be dealt with similarly with the fact that $S_{\mu,\,\nu}(\zeta)$ is an even function of $\nu$.
\smallskip

\begin{lemma}\label{L:lommel-continuation-general-1}
If $-\nu \not\in\{0,\,1,\,2,\ldots\}$, then for each integer $m$,

\begin{align}\label{E:lommel-continuation-general-2}
S_{\nu-1,\,\nu}(\zeta {\rm e}^{-m\pi i})&={\rm e}^{-m\nu\pi
i}S_{\nu-1,\,\nu}(\zeta)+K'_+ H_\nu^{(1)}(\zeta)+K'_-
H_\nu^{(2)}(\zeta),
\end{align}\smallskip

\noindent where

\begin{equation}\label{E:K-prime}
K'_\pm=\pi 2^{\nu-2} i{\rm e}^{-m\nu\pi i}\Gamma(\nu)\big[U_{m-1}(\cos \nu\pi){\rm e}^{(m \pm 1)\nu\pi i}-m\big].
\end{equation}
\end{lemma}

\begin{proof} Since $-\nu \not\in\{0,\,1,\,2,\ldots\}$, it follows from the formulae (\ref{E:lommel-(nu-1)}),
(\ref{E:bessel-first-continuation}--\ref{E:bessel-second-continuation})
that

\begin{align*}
S_{\nu-1,\,\nu}(\zeta {\rm e}^{-m\pi
i})&=2^{\nu-1}\Gamma(\nu)J_\nu(\zeta {\rm e}^{-m\pi i})\log
(\zeta {\rm e}^{-m\pi i})-2^{\nu-2}\pi \Gamma(\nu)Y_\nu(\zeta {\rm
e}^{-m\pi i}) \notag\\
&\quad-\frac{1}{4}(\zeta {\rm e}^{-m\pi i})^\nu\Gamma(\nu)\sum_{k=0}^{+\infty}\frac{(-1)^k
(\frac{1}{2}\zeta {\rm e}^{-m\pi i})^{2k}A(k)}{k!\Gamma(\nu+k+1)}\notag\\
&=2^{\nu-1} \Gamma(\nu)(\log \zeta-m\pi i){\rm e}^{-m\nu\pi
i}J_\nu(\zeta)\\
&\quad-2^{\nu-2}\pi \Gamma(\nu)[{\rm e}^{m\nu\pi
i}Y_\nu(\zeta)-2i\sin(m\nu\pi)\cot(\nu\pi)J_\nu(\zeta)]\\
&\quad+{\rm e}^{-m\nu\pi
i}\big[S_{\nu-1,\,\nu}(\zeta)-2^{\nu-1}\Gamma(\nu)J_\nu(\zeta)\log\zeta
+2^{\nu-2}\pi\Gamma(\nu)Y_\nu(\zeta)\big]\\
&={\rm e}^{-m\nu\pi i}S_{\nu-1,\,\nu}(\zeta)+2^{\nu-2}\pi
\Gamma(\nu){\rm e}^{-m\nu\pi i}Y_\nu(\zeta)\\
&\quad -2^{\nu-1}m\pi i\Gamma(\nu){\rm e}^{-m\nu\pi i}J_\nu(\zeta)\\
&\quad-2^{\nu-2}\pi \Gamma(\nu)[{\rm e}^{m\nu\pi i}Y_\nu(\zeta)-2i\sin(m\nu\pi)\cot(\nu\pi)J_\nu(\zeta)]\\
 \medskip &={\rm e}^{-m\nu\pi i}S_{\nu-1,\,\nu}(\zeta)-2^{\nu-1}\pi i\Gamma(\nu)\sin
(m\nu\pi) Y_\nu(\zeta)\\
&\quad+2^{\nu-1}\pi i \Gamma(\nu)\big[\sin(m\nu\pi)\cot(\nu\pi)-m{\rm e}^{-m\nu\pi i}\big]J_\nu(\zeta)\\
\medskip
&={\rm e}^{-m\nu\pi i}S_{\nu-1,\,\nu}(\zeta)+K'_+
H_\nu^{(1)}(\zeta)+K'_- H_\nu^{(2)}(\zeta),
\end{align*}

\noindent where $m$ is an integer and

\begin{align*}
K'_\pm=2^{\nu-2}\pi \Gamma(\nu) \big\{\mp \sin(m\nu \pi)+i\big[\sin(m \nu \pi)\cot(\nu \pi)-m{\rm e}^{-m\nu\pi i}\big]\big\}.
\end{align*}\smallskip

It is a routine verification that the above expression for $K'_{\pm}$ can be reduced to our desired form (\ref{E:K-prime}), and thus proving the formula (\ref{E:K-prime}).
\end{proof}

It remains to substitute the above formula for $S_{\nu-1,\,\nu}(\zeta)$ into (\ref{E:lommel-negative}) to obtain a continuation formula of $S_{\nu-2p-1,\,\nu}(\zeta)$ when $-\nu \not\in\{0,\,1,\,2,\ldots\}$:

\begin{lemma}\label{L:lommel-continuation-general-2}
If $-\nu\not \in \{0,\,1,\,2,\ldots\}$, then for each integer $m$,

\begin{align*}
S_{\nu-2p-1,\,\nu}(\zeta {\rm e}^{-m\pi i})&={\rm e}^{-m\nu\pi
i}S_{\nu-2p-1,\,\nu}(\zeta)+\frac{(-1)^pK_+'}{2^{2p}p!(1-\nu)_p}H_\nu^{(1)}(\zeta)\\
&\quad+ \frac{(-1)^pK_-'}{2^{2p}p!(1-\nu)_p}H_\nu^{(2)}(\zeta).
\end{align*}
\end{lemma}

The above consideration involving the { Lemmae \ref{L:lommel-continuation-general-1} and \ref{L:lommel-continuation-general-2}} which deal with analytic continuation formulae when $-\nu \not\in\{0,\,1,\,2,\ldots\}$. We next treat the {case that} $\nu=0$ in the following results.

\begin{lemma}\label{L:lommel-continuation-general-3} If $\nu=0$, then for any integer $m$,

\begin{align}\label{E:lommel-continuation-general-3}
S_{-1,\,0}(\zeta {\rm e}^{-m\pi
i})=S_{-1,\,0}(\zeta)+K_+''H_0^{(1)}(\zeta)+K_-''H_0^{(2)}(\zeta),
\end{align}\smallskip

\noindent where

$$K_\pm''=-\frac{m\pi^2(m \pm 1)}{4}.$$
\end{lemma}

\begin{proof} We have $Y_0(\zeta {\rm e}^{-m\pi i})=Y_0(\zeta)-2mi
J_0(\zeta)$ from the equation (\ref{E:bessel-second-continuation}),
so (\ref{E:lommel-negative-one-further}) gives

\begin{align*}
S_{-1,\,0}(\zeta {\rm e}^{-m\pi i})
&=-\frac{1}{2}J_0(\zeta)(\log \zeta -m\pi
i)^2+\frac{\pi}{2}[Y_0(\zeta)-2mi J_0(\zeta)](\log \zeta -m\pi i)\\
&\quad+\frac{1}{2}\sum_{k=0}^{+\infty}\frac{(-1)^k(\frac{1}{2}\zeta)^{2k}B(k)}{(k!)^2}\notag\\
&=S_{-1,\,0}(\zeta)-\frac{m^2\pi^2}{2}J_0(\zeta)-\frac{m\pi^2}{2}iY_0(\zeta).
\end{align*}\smallskip

\noindent Hence the formula (\ref{E:lommel-continuation-general-3}) follows
from (\ref{E:hankels-definition}).
\end{proof}

\begin{lemma}\label{L:lommel-continuation-general-4}
If $\nu=0$, then for any integer $m$,

\begin{equation*}
S_{-2p-1,\,0}(\zeta {\rm e}^{-m\pi
i})=S_{-2p-1,\,0}(\zeta)+\frac{(-1)^pK_+''}{2^{2p}(p!)^2}H_0^{(1)}(\zeta)+
\frac{(-1)^pK_-''}{2^{2p}(p!)^2}H_0^{(2)}(\zeta).
\end{equation*}
\end{lemma}
\medskip

\begin{proof} This is easily obtained by substituting $\zeta {\rm e}^{-m\pi
i}$ into the equation (\ref{E:lommel-negative}) with $\nu=0$ and applying
(\ref{E:lommel-continuation-general-3}).
\end{proof}\smallskip

Finally, the case when $-\nu \in \mathbf{N}$ is now considered in the next two results.

\begin{lemma}\label{L:lommel-continuation-odd-negative} Suppose that $m$ is any integer. We define $\delta_m=1+(-1)^{m-1}$ and for every polynomial $P_n(\zeta)$ of degree $n$, we define $\widehat{P}_n(\zeta)$ to be the polynomial containing the term of $P_n(\zeta)$ with odd powers in $\zeta$ and $\overline{P}_n(\zeta):=P_n(\zeta)-\delta_m\widehat{P}_n(\zeta)$. If $\nu=-n$ for a positive integer $n$, then we have


\begin{align}\label{E:lommel-continuation-odd-negative1}
S_{-n-1,\,-n}(\zeta {\rm e}^{-m\pi i})&=(-1)^{mn}S_{-n-1,\,-n}(\zeta)+\frac{(-1)^{(m+1)n}}{2^n n!}\zeta^{-n}\times\\
&\qquad\Big\{-\delta_m\big[\widehat{A}_n(\zeta)+\widehat{B}_n(\zeta)S_{-1,\,0}(\zeta)+\zeta\widehat{C}_n(\zeta)S_{-1,\,0}'(\zeta)\big] \notag\\
&\qquad+\overline{B}_n(\zeta)\big[K_+''H_0^{(1)}(\zeta)+K_-''H_0^{(2)}(\zeta)\big] \notag\\ &\qquad-\zeta\overline{C}_n(\zeta)\big[K_+''H_1^{(1)}(\zeta)+K_-''H_1^{(2)}(\zeta)\big]\Big\}.\notag
\end{align}
\end{lemma}\smallskip

\begin{proof} Differentiating (\ref{E:lommel-continuation-general-3}) and applying (\ref{E:hankel-derivatives}) yields


\begin{equation}\label{E:lommel-continuation-odd-negative2}
{\frac{\rm d}{\rm d\zeta}}S_{-1,\,0}(\zeta {\rm e}^{-m\pi i})=S^\prime_{-1,\,0}(\zeta)-K_+''H_1^{(1)}(\zeta)-K_-''H_1^{(2)}(\zeta).
\end{equation}\smallskip

We note that it can be derived from the definition easily that for every polynomial $P_n(\zeta)$ of degree $n$ and every integer $m$, we must have

\begin{equation}\label{E:lommel-continuation-odd-negative3}
P_n(\zeta {\rm e}^{-m\pi i})=P_n(\zeta)-\delta_m\widehat{P}_n(\zeta)=\overline{P}_n(\zeta).
\end{equation}
\smallskip

\noindent Hence it can be seen without difficulty that the formula (\ref{E:lommel-continuation-odd-negative1})
follows from the formulae (\ref{E:lommel-negative-integer-general}), { (\ref{E:lommel-continuation-general-3})}, (\ref{E:lommel-continuation-odd-negative2}) and (\ref{E:lommel-continuation-odd-negative3}).

\end{proof}\smallskip

Now we can substitute the formula (\ref{E:lommel-continuation-odd-negative1}) into the formula (\ref{E:lommel-negative}) to obtain a continuation formula of $S_{-n-2p-1,\,-n}(\zeta)$:

\begin{lemma}\label{L:lommel-continuation-odd-negative-general} If $\nu=-n$ for a positive integer $n$, then for any integer $m$,

\begin{align}\label{E:lommel-continuation-odd-negative-general}
S_{-n-2p-1,\,-n}&(\zeta {\rm e}^{-m\pi i})\\
&=(-1)^{mn}S_{-n-2p-1,\,-n}(\zeta)+\frac{(-1)^{(m+1)n+p}}{2^{2p+n} n!p!(1+n)_p}\zeta^{-n}\times \notag\\
&\qquad\Big\{-\delta_m\big[\widehat{A}_n(\zeta)+\widehat{B}_n(\zeta)S_{-1,\,0}(\zeta)+\zeta\widehat{C}_n(\zeta)S_{-1,\,0}'(\zeta)\big] \notag\\
&\qquad+\overline{B}_n(\zeta)\big[K_+''H_0^{(1)}(\zeta)+K_-''H_0^{(2)}(\zeta)\big]\notag\\
&\qquad-\zeta\overline{C}_n(\zeta)\big[K_+''H_1^{(1)}(\zeta)+K_-''H_1^{(2)}(\zeta)\big]\Big\}.\notag
\end{align}
\end{lemma}\medskip

\subsection{An asymptotic expansion of $S_{\mu,\,\nu}(\zeta)$}
It is known that when $\mu \pm \nu$ are not odd positive integers,
then $S_{\mu,\,\nu}(\zeta)$ has the asymptotic expansion

\begin{align}\label{E:lommel-asymptotic-expansion}
    S_{\mu,\,\nu}(\zeta)
    =\zeta^{\mu-1}\left[\sum_{k=0}^p \frac{(-1)^k
    c_k}{\zeta^{2k}}\right] +O\left(\zeta^{\mu-2p-2}\right)
\end{align}\smallskip

\noindent for large $|\zeta|$ and $|\arg \zeta|<\pi$, where $p$ is a
non-negative integer and the numbers $c_k$ are the coefficients
defined in (\ref{ck}). See also \cite{Wa44}, \S10.75.

\begin{remark}
\label{R:lommel-terminate}\rm It is clear that
(\ref{E:lommel-asymptotic-expansion}) is a series in descending
powers of $\zeta$ starting from the term $\zeta^{\mu-1}$ and
(\ref{E:lommel-asymptotic-expansion}) terminates if one of the
numbers $\mu \pm \nu$ is an odd positive integer. In particular,
if $\mu-\nu=2p+1$ for some non-negative integer $p$, then we have $K_+=0$ in the analytic continuation formula (\ref{E:lommel-continuation-general}) and thus, in this degenerate case, the formula (\ref{E:lommel-continuation-general}) becomes

\begin{equation*}
S_{2p+1+\nu,\,\nu}(\zeta {\rm e}^{-m\pi i})={\rm e}^{-m\nu\pi i} S_{2p+1+\nu,\,\nu}(\zeta)
\end{equation*}\smallskip

\noindent for every integer $m$ and $|\arg \zeta|<\pi$.
\end{remark}
\medskip

\subsection{Linear independence of Lommel's functions}
We next discuss the linear independence of the Lommel
functions $S_{\mu_j,\,\nu}(\zeta)$.\par

\begin{lemma}\label{E:lommel-independence-1} Suppose $n \ge 2$, and $\mu_j$
and $\nu$ be complex numbers such that $\Re(\mu_j)$ are
all distinct for $j=1,2,\ldots, n$. Then the Lommel functions

$$S_{\mu_1,\, \nu}(\zeta),\,S_{\mu_2,\,
\nu}(\zeta),\ldots,\,S_{\mu_n,\, \nu}(\zeta)$$\smallskip

\noindent are linearly independent.
\end{lemma}

\begin{proof}
Let us now assume that the Lommel functions $S_{\mu_j,\,\nu}(\zeta), j=1,2,\ldots, n$, to be linearly dependent.
Then there exist constants $C_j$ not all zero such that

\begin{equation}
\label{E:lommel-dependence} \D \sum_{j=1}^{n}
C_jS_{\mu_j,\,\nu}(\zeta)=0.
\end{equation}\smallskip

\noindent We may assume, without loss of generality, that none of the constants $C_j$ is zero, and  that
$\Re(\mu_1)<\Re(\mu_2)<\cdots<\Re(\mu_{n})$.\smallskip

We substitute the asymptotic expansions (\ref{E:lommel-asymptotic-expansion}) of the Lommel functions into
(\ref{E:lommel-dependence}) and consider only the leading terms there. We deduce

\[
\D-C_{n}S_{\mu_{n},\, \nu}(\zeta)=\sum_{j=1}^{n-1}C_jS_{\mu_j,\, \nu}(\zeta),
\]\smallskip

\noindent and

\begin{align*}
\D |C_{n}||\zeta|^{\Re(\mu_{n})-1} {\rm e}^{-\Im(\mu_{n})\arg \zeta}
&\le
\sum_{j=1}^{n-1}|C_j||\zeta|^{\Re(\mu_j)-1}{\rm e}^{-\Im(\mu_j)\arg \zeta},\\
|C_{n}||\zeta|^{\Re(\mu_{n})-1}{\rm e}^{-\Im(\mu_{n})\arg
\zeta} &\le (n-1)|C_{n-1}||\zeta|^{\Re(\mu_{n-1})-1}{\rm e}^{-\Im(\mu)\arg \zeta},\\
|C_{n}| {\rm e}^{-\Im(\mu_{n})\arg \zeta}&\le (n-1)
|C_{n-1}||\zeta|^{[\Re(\mu_{n-1})-\Re(\mu_{n})]}{\rm
e}^{-\Im(\mu)\arg \zeta}
\end{align*}\smallskip

\noindent for sufficiently large $\zeta$ and $|\arg \zeta|<\pi$,
where
$\Im(\mu)=\min\big\{\Im(\mu_1),\,\ldots,\,\Im(\mu_{n-1})\big\}$.
Since $\Re(\mu_{n-1})-\Re(\mu_{n})<0$, the right hand side of the
last inequality approaches zero as $\zeta \to \infty$ in $|\arg
\zeta|<\pi$ which is a contradiction. Hence, our desired result
follows.
\end{proof}
\medskip

\section{\bf Proof of Theorem {\rm \ref{T:Chiang-Yu}(b)}}
\label{P:Chiang-Yu}
The general form of the solution of (\ref{E:Chiang-Yu}) was already
derived in \S 2.\smallskip

Let $y(\zeta)$ be the general solution of
(\ref{E:generalized-Lommel}). We shall recall from Theorem \ref{T:Chiang-Yu} that $A,\,B,\, L,\, M,\,\sigma_j$ and
$\mu_j$ are complex constants, $L$ and $M$ are non-zero and at least one of $\sigma_j,\,j=1,\,2,\ldots,\,n$, being non-zero. Thus it follows from the definitions
(\ref{E:hankels-definition}) that we may rewrite the general
solution (\ref{E:general-soln-bessel}) of the equation
(\ref{E:generalized-Lommel}) in terms of Hankel's functions in the
form

\begin{equation} \label{E:general-soln-hankel}
y(\zeta)=C
H_\nu^{(1)}(\zeta)+DH_\nu^{(2)}(\zeta)+\sum_{j=1}^n \sigma_j
S_{\mu_j,\, \nu}(\zeta),
\end{equation}\smallskip

\noindent where $C=\frac{1}{2}(A-iB)$ and $D=\frac{1}{2}(A+iB)$. It is easily seen that $A=B=0$ if and only if $C=D=0$. Thus the general solution $f(z)={\rm e}^{-N z}y(L {\rm e}^{M z})$ of (\ref{E:Chiang-Yu}) assumes the form

\begin{align}\label{E:f-solution}
f(z)={\rm e}^{-Nz}\bigg[CH_\nu^{(1)}(L{\rm e}^{Mz})+D{\rm e}^{-Nz}H_\nu^{(2)}(L{\rm e}^{Mz})+\sum_{j=1}^n \sigma_jS_{\mu_j,\, \nu}(L{\rm e}^{Mz})\bigg].
\end{align}\smallskip

To prove Theorem \ref{T:Chiang-Yu}(b), we must \textit{first} show that if $f(z)$ is subnormal, then we have $C=D=0$. The proof of this depends on the following result:

\begin{proposition}\label{T:G(z)-growth}
Suppose $C$ and $D$ are complex numbers such that $(C,\, D)\not=(0,\,0)$. Then there exists a sequence of complex numbers $\{z_n\}=\big\{r_n{\rm e}^{i\theta_n}\big\}$ such that $|z_n|=r_n\to
+\infty$ as $n\to+\infty$, $-\pi<\arg (L{\rm e}^{Mz_n})<\pi$ for all positive integers $n$ and that the entire function

\begin{align}\label{E:G(z)}
G(z)&= CH_\nu^{(1)}(L{\rm e}^{Mz})+DH_\nu^{(2)}(L{\rm e}^{Mz})+\sum_{j=1}^n \sigma_jS_{\mu_j,\, \nu}(L{\rm e}^{Mz})\notag\\
&:= F(z)+\sum_{j=1}^n \sigma_jS_{\mu_j,\, \nu}(L{\rm e}^{Mz})
\end{align}\smallskip

\noindent satisfies the estimate
	\begin{align}\label{E:maximum-modulus}
		M(r_n,\,G) & {\ge} \big({2^{1/2}}{{\pi}^{-1}}R_n^{-1}\big)^\frac12\notag\\
	&\qquad\times\left\{
		\begin{array}{ll}
		|C{\rm e}^{-i\frac{\nu\pi}{2}}|{\rm e}^{R_n} { \bigg(\D 1-\bigg|\frac{D}{C}{\rm e}^{i\nu\pi}\bigg|{\rm e}^{-2R_n}\bigg)}+o({\rm e}^{R_n}),
		& \text{if\ } C\not=0;\\
        & \\
		|D {\rm e}^{i\frac{\nu\pi}{2}}|{\rm e}^{R_n}+o({\rm e}^{R_n}),
		& \text{if\ } C=0;\\
		\end{array}
		\right.
	\end{align}

\noindent where $\D R_n=\frac{|L|}{\sqrt{2}}{\rm e}^{\frac{|M|r_n}{\sqrt{2}}}$. In particular, the entire function {\rm (\ref{E:f-solution})} is not subnormal.
\end{proposition}
\medskip

\begin{proof}[Proof of Proposition \ref{T:G(z)-growth}] We shall estimate the growth of each of the functions in the (\ref{E:G(z)}) in the two lemmae below.

\begin{lemma}\label{L:F-not-subnormal}
Suppose $C$ and $D$ are complex numbers such that $(C,\, D)\not=(0,\,0)$. Then there
exists a sequence $\{z_n\}=\big\{r_n{\rm e}^{i\theta_n}\big\}$ such that
$r_n \to +\infty$ as $n \to +\infty$, $-\pi <\arg (L{\rm e}^{Mz_n})<\pi$, where $\theta_n$ is fixed for all positive integers $n$, and that the entire function

\begin{equation}\label{E:lemma5.1}
F(z)=CH_\nu^{(1)}(L {\rm e}^{M z})+DH_\nu^{(2)}(L {\rm e}^{M z})
\end{equation} \smallskip

\noindent satisfies the estimate {\rm (\ref{E:maximum-modulus})} with $G$ replaced by $F$.
Hence we have $\sigma(F)=+\infty$ and the $F$ is not subnormal from the
definition {\rm (\ref{E:subnormal})}.
\end{lemma}

\begin{proof}[Proof of Lemma \ref{L:F-not-subnormal}] We consider the growth
of (\ref{E:lemma5.1}) in the \textit{principal branch} of
$H_\nu^{(1)}(\zeta)$ and $H_\nu^{(2)}(\zeta)$.\footnote{The
principal branch of $H_\nu^{(1)}(\zeta)$ and $H_\nu^{(2)}(\zeta)$ is
assumed to be $-\pi<\arg \zeta<\pi$.} Let

$$L=|L|{\rm e}^{ia},\quad M=|M|{\rm e}^{ib} \quad \mbox{and}\quad
z_n=r_n{\rm e}^{i\theta_n},$$\smallskip

\noindent where $a,\, b \in (-\pi,\pi {]}$ and $r_n \to +\infty$ as $n \to +\infty$. The idea of proof is to choose
suitable sequences $\{\theta_n\}$ and $\{r_n\}$ (and hence $\{z_n\}$), so that we can apply the asymptotic expansions
(\ref{E:hankel1-asy}) and (\ref{E:hankel2-asy}) {\it simultaneously} to estimate explicitly the growths of $H_\nu^{(1)}(L{\rm e}^{Mz_n})$ and $H_\nu^{(2)}(L{\rm e}^{Mz_n})$.\smallskip

We consider the sequence

\begin{equation}\label{E:sequence-z}
\{z_n\}=\big\{r_n {\rm e}^{i {\theta_n}}\big\},
\end{equation}

\noindent where $\D \theta_n=\frac{\pi}{4}-b$ for all $n \in
\mathbf{N}$ and

\begin{equation}\label{E:sequence-r_n}
r_n=\left\{
  \begin{array}{ll}
    \D\frac{\sqrt{2}}{|M|}\Big(2n\pi-\frac{\pi}{4}-a\Big), & \hbox{if $C\not=0$;} \\
    \D\frac{\sqrt{2}}{|M|}\Big(2n\pi+\frac{\pi}{4}-a\Big), & \hbox{if $C=0$.}
  \end{array}
\right.
\end{equation}\smallskip

\noindent Hence a routine computation yields

\begin{align}\label{E:sequence-e^z}
L{\rm e}^{Mz_n}&=|L|\exp\big\{|M|r_n\cos(b+\theta_n)+i\big[a+|M|
r_n\sin(b+\theta_n)\big]\big\} \notag\\
&=|L|\exp\bigg[\frac{|M|r_n}{\sqrt{2}}+i\bigg(a+\frac{|M|r_n}{\sqrt{2}}\bigg)\bigg] \notag\\
&=|L|{\rm e}^{\frac{|M|r_n}{\sqrt{2}}}{\rm e}^{i(2n\pi \mp \frac{\pi}{4})} \notag\\
&:=\left\{
  \begin{array}{ll}
    \D R_n(1-i), & \hbox{if $C\not=0$;} \\
    \D R_n(1+i), & \hbox{if $C=0$;}
  \end{array}
\right.
\end{align}\smallskip

\noindent {and} it is easy to see that $-\pi<\arg(L{\rm e}^{Mz_n})<\pi$ for each positive integer $n$. Now we can apply the asymptotic expansions (\ref{E:hankel1-asy}) and (\ref{E:hankel2-asy}) for sufficiently large $n$ and (\ref{E:sequence-e^z}) to obtain the following estimates:

\begin{align}
	F(z_n) &=CH_\nu^{(1)}(L{\rm e}^{Mz_n})+D H_\nu^{(2)}(L{\rm e}^{Mz_n}) \notag\\
	&=\left(\frac{2}{\pi R_n(1\mp i)}\right)^{1/2}\left\{
	{ \frac{C(1-i)}{\sqrt{2}}}{\rm e}^{-i\frac{\nu\pi}{2}}{\rm e}^{iR_n(1\mp i)}
	\big(1+O(R_n^{-1})\big)\right.\notag\\
	&\qquad\quad\left.+{ \frac{D(1+i)}{\sqrt{2}}} {\rm e}^{i\frac{\nu\pi}{2}}{\rm e}^{-iR_n(1\mp i)}	 \big(1+O(R_n^{-1})\big)\right\}\notag\\
	&=\left\{
		\begin{array}{ll}
			\D { \left(\frac{2}{\pi R_n(1 - i)}\right)^{1/2}\frac{C(1-i)}{\sqrt{2}}}{\rm
            e}^{-i\frac{\nu\pi}{2}}{\rm e}^{iR_n(1- i)}&\\
			\qquad\D \times \Bigg[1+\frac{D(1+i)}{C(1-i)}{\rm e}^{i\nu\pi} {\rm e}^{-2iR_n(1- i)}\Bigg]+o\big({ {\rm e}^{R_n+iR_n}}\big)&\\
            \qquad +o\big({\rm e}^{-R_n-iR_n}\big), &\text{if\ } C\not=0;\\
            & \\
			\D { \left(\frac{2}{\pi R_n(1+ i)}\right)^{1/2}\frac{D(1+i)}{\sqrt{2}}} {\rm e}^{i\frac{\nu\pi}{2}}{\rm e}^{-iR_n(1+ i)}+o\big({\rm e}^{R_n-iR_n}\big),
			& \text{if\ } C=0.
		\end{array}
		\right.
\end{align}

\noindent Hence

	\begin{align}\label{E:estimate-hankel12}
		M(r_n,\,F) & {\ge} |F(z_n)|\notag\\
		&{\ge} \big({2^{1/2}}{{\pi}^{-1}}R_n^{-1}\big)^\frac12\notag\\
	&\qquad\times\left\{
		\begin{array}{ll}
		\D |C{\rm e}^{-i\frac{\nu\pi}{2}}|{\rm e}^{R_n} {\bigg(1-\bigg|\frac{D}{C} {\rm e}^{i\nu\pi}\bigg|{\rm e}^{-2R_n}\bigg)}+o({\rm e}^{R_n}),
		& \text{if\ } C\not=0;\\
& \\
		|D{\rm e}^{i\frac{\nu\pi}{2}}|{\rm e}^{R_n}+o({\rm e}^{R_n}),
		& \text{if\ } C=0.\\
		\end{array}
		\right.
	\end{align}

\noindent
Hence the estimate (\ref{E:estimate-hankel12}) is our desired result. Clearly the (\ref{E:estimate-hankel12}) also implies that $\sigma(F)=+\infty$.
\smallskip

The same estimate (\ref{E:estimate-hankel12}) also shows that, as $n \to +\infty$

$$\frac{\log\log{\rm e}^{R_n}}{ r_n}=\frac{\D
\log\frac{|L|}{\sqrt{2}}+\frac{|M|r_n}{\sqrt{2}}}{r_n} \to
\frac{|M|}{\sqrt{2}}\not=0.$$\smallskip

\noindent It follows from the definition (\ref{E:subnormal})  that
$F(z)$ is not subnormal. This completes the proof of Lemma \ref{L:F-not-subnormal}.
\end{proof}
\medskip

We next estimate the growth of the Lommel function $S_{\mu,\,
\nu}(L{\rm e}^{Mz})$ on the same sequence $\{z_n\}$ defined in Lemma \ref{L:F-not-subnormal}.
\medskip

\begin{lemma}\label{L:lommel-modulus}
Let $S_{\mu,\, \nu}(\cdot)$ be the principal branch of the Lommel function, where $\mu \pm \nu$ are not odd positive integers, { and} $L,\, M$ are {non-zero} constants. Then on the sequence {\rm (\ref{E:sequence-z})} defined in Lemma {\rm 4.2}, we have

\begin{align*}
\big|S_{\mu,\,\nu}(L{\rm e}^{Mz_n})\big| \sim \big|\big(L{\rm
e}^{Mz_n}\big)^{\mu-1}\big|={ \big(|L|{\rm e}^{\frac{|M|r_n}{\sqrt{2}}}\big)^{[\Re (\mu)-1]}} \times {\rm
e}^{\varepsilon \frac{\pi}{4}\Im(\mu)},
\end{align*}\smallskip

\noindent where the value of $\varepsilon=\pm 1$ depends on the sequence {\rm (\ref{E:sequence-r_n})} chosen.
\end{lemma}
\medskip

\begin{proof}[Proof of Lemma \ref{L:lommel-modulus}] It is clear that $-\pi<\arg(L{\rm e}^{Mz_n})<\pi$ for the sequence (\ref{E:sequence-z}) (see also (\ref{E:sequence-e^z})). Thus, by choosing $p=0$ in the asymptotic
expansion (\ref{E:lommel-asymptotic-expansion}), we obtain that

\begin{align*}
|S_{\mu,\,\nu}(L{\rm e}^{Mz_n})| &\sim \big|\big(L{\rm e}^{Mz_n}\big)^{\mu-1}\big|\\
&=\big|\big(L{\rm e}^{Mz_n}\big)^{\Re(\mu)-1}\big|\cdot \big|\big(L{\rm e}^{Mz_n}\big)^{i\Im(\mu)}\big|\\
&={ \big(|L|{\rm e}^{\frac{|M|r_n}{\sqrt{2}}}\big)^{[\Re (\mu)-1]}} \cdot\big|\big(L{\rm e}^{Mz_n}\big)^{i\Im(\mu)}\big|.
\end{align*}\smallskip

\noindent It is clear from (\ref{E:sequence-e^z}) that

\begin{align*}
\bigg|\big(L{\rm e}^{Mz_n}\big)^{i{\Im(\mu)}}\bigg|=\big|{\rm e}^{i \Im(\mu)\log [R_n(1 \mp i)]}\big|=\big|{\rm e}^{-\Im(\mu)\arg (1 \mp i)}\big|={\rm e}^{\varepsilon \frac{\pi}{4}{\Im(\mu)}},
\end{align*}\smallskip

\noindent where the value of $\varepsilon=\pm 1$ depends on the sequence (\ref{E:sequence-r_n}) such that $\varepsilon=+1$ if $C\not=0$ and $\varepsilon=-1$ otherwise.
Thus we complete the proof of Lemma \ref{L:lommel-modulus}.
\end{proof}

We can now prove Proposition \ref{T:G(z)-growth}. Let the right hand side of
(\ref{E:G(z)}) be in the principal branch of
$H_\nu^{(1)}(\zeta),\,H_\nu^{(2)}(\zeta)$ and
$S_{\mu_j,\,\nu}(\zeta)$ for $j=1,2,\ldots,n$.\footnote{The
principal branch of each $S_{\mu_j,\,\nu}(\zeta)$, $j=1,2,\ldots,n$,
is assumed to be $-\pi<\arg \zeta<\pi$.} It is obvious that

$$|G(z)|\ge \big|F(z)\big|-\sum_{j=1}^n \big|\sigma_j S_{\mu_j,\,
\nu}(L{\rm e}^{Mz})\big|.$$\smallskip

\noindent If in addition that we let $z=z_n$ be the sequence (\ref{E:sequence-z}), then the estimates in Lemmae \ref{L:F-not-subnormal} and \ref{L:lommel-modulus} clearly imply that $G(z)$ satisfies the estimate (\ref{E:maximum-modulus}) for all sufficiently large $n$ in the
principal branch of the functions
$H_\nu^{(1)}(\zeta),\,H_\nu^{(2)}(\zeta)$ and
$S_{\mu_j,\,\nu}(\zeta)$. \smallskip

By the similar argument as in the proof of Lemma \ref{L:F-not-subnormal}, we know that
$G(z)$ is not subnormal. Since we have

$$f(z)={\rm e}^{-Nz}\Bigg[F(z)+\sum_{j=1}^n
\sigma_jS_{\mu_j,\,\nu}(L{\rm e}^{Mz})\Bigg]={\rm
e}^{-Nz}G(z)$$\smallskip

\noindent and the function ${\rm e}^{-Nz}$ is clearly subnormal,
the function $f(z)$ is not subnormal too.
\end{proof}\smallskip

We now continue the proof of Theorem \ref{T:Chiang-Yu}(b). Since $f$ is subnormal, so according to the above analysis we must have $C=D=0$. That is, $f(z)$ must reduce to the following form:

\begin{equation}\label{E:lommel-sum}
f(z)=\sum_{j=1}^n \sigma_j {\rm e}^{-N z}S_{\mu_j,\, \nu}(L {\rm
e}^{M z}).
\end{equation}\smallskip

\noindent Since Lemma \ref{E:lommel-independence-1} shows that the Lommel functions $S_{\mu_j,\,\nu} (\zeta), \, j=1,2, \ldots, n$, are linearly independent over $\mathbf{C}$ and not all $\sigma_i$ are zero, the solution (\ref{E:lommel-sum}) is clearly not identically zero.
\vspace{0.5cm}

To complete the proof of Theorem \ref{T:Chiang-Yu}(b), we also need to prove that
when $\sigma_j$ is non-zero, $\mu_j$
and $\nu$ must satisfy either

\begin{equation}\label{E:desired-equations}
\cos \bigl(\frac{\mu_j+\nu}{2}\pi\bigr)=0\quad \mbox{or}\quad 1+{\rm
e}^{-(\mu_j+\nu)\pi i}=0,
\end{equation}\smallskip

\noindent where $j \in \{1,\, 2,\ldots,\, n\}$. To do so we will
need Lemma \ref{L:lommel-modulus} and the following result.

\begin{proposition}\label{L:lommel-not-subnormal}
Let $p$ be a non-negative integer and let $\nu$ be an arbitrary complex number such that if $\nu$ is { an integer}, then it is not greater than $p$. Then the entire function $S_{\nu-2p-1,\,\nu}(L{\rm e}^{Mz})$ is not subnormal.\footnote{By the first paragraph of \S 3.3, the Proposition \ref{L:lommel-not-subnormal} is also valid when $\mu+\nu$ is an odd negative integer.}
\end{proposition}

\begin{proof}[Proof of Proposition \ref{L:lommel-not-subnormal}] We assume that
$S_{\nu-2p-1,\,\nu}(L{\rm e}^{Mz})$ is subnormal. We recall from
the beginning of \S3 that its growth must be independent of the
different branches of $S_{\nu-2p-1,\,\nu}(\zeta)$. Let us distinguish three cases:

	\begin{enumerate}
		\item[(i)] $-\nu\not\in \{0,\,1,\,2,\ldots\}$. It follows from this branches-argument and the Lemma \ref{L:lommel-continuation-general-2} that

        \begin{align}\label{E:lommel-nu-non-zero}
        S_{\nu-2p-1,\,\nu}(L{\rm e}^{Mz}{\rm e}^{-m\pi i})&={\rm e}^{-m\nu\pi i}S_{\nu-2p-1,\,\nu}(L{\rm e}^{Mz}) \\
        &\quad+\frac{(-1)^pK_+'}{2^{2p}p!(1-\nu)_p}H_\nu^{(1)}(L{\rm e}^{Mz}) \notag \\
        &\quad+\frac{(-1)^pK_-'}{2^{2p}p!(1-\nu)_p}H_\nu^{(2)}(L{\rm
        e}^{Mz}),\notag
        \end{align}
where ${K'_\pm}$ are given by (\ref{E:K-prime}),
 is required to be subnormal for \textit{each} integer $m$. We note that the Lommel and Hankel functions on the right side of (\ref{E:lommel-nu-non-zero}) are in their principal branch. Let $G(z):=S_{\nu-2p-1,\,\nu}(L{\rm e}^{Mz}{\rm e}^{-m\pi i})$ in {the} Proposition \ref{T:G(z)-growth}. Then we deduce that this $G(z)$ also satisfies the estimate (\ref{E:maximum-modulus}) on the sequence (\ref{E:sequence-z}). Thus it is not subnormal unless $\D \frac{K'_\pm}{(1-\nu)_p}=0$ and  from which we deduce

        \begin{equation}\label{E:K'_pm-non-zero}
        \frac{\sin (m\nu\pi)}{(1-\nu)_p}=0
        \end{equation}\smallskip

        \noindent for \textit{each} integer $m$. If $\nu=n \in \{1,\,2,\,\ldots,\,p\}$, then {equation (\ref{E:K'_pm-non-zero})} implies

        \begin{equation*}
        0=\lim_{\nu \to n}\frac{\sin (m\nu\pi)}{(1-\nu)_p}=\frac{m\pi(-1)^{mn}}{-(1-n)(2-n)\cdots(-1)(1)\cdots(p-n)}
        \end{equation*}\smallskip

        \noindent which is valid only when $m=0$, a contradiction. Therefore we must have $\nu \not\in\{1,\,2,\,\ldots,\,p\}$ and equation (\ref{E:K'_pm-non-zero}) implies that $\sin (m\nu\pi)=0$ for every integer $m$. Thus $\nu$ is an integer and $\nu=n \ge p+1$, a contradiction to the assumption. Hence we have $\D \frac{K'_\pm}{(1-\nu)_p} \neq 0$ and so $S_{\nu-2p-1,\,\nu}(L{\rm e}^{Mz})$ is not subnormal in this case by the Proposition \ref{T:G(z)-growth}.\smallskip

        We now consider the second case.\smallskip

    \item[(ii)] If $\nu=n=0$, then the independence of the branches means that we need to apply the analytic continuation formula in Lemma \ref{L:lommel-continuation-general-4}  in our consideration instead. Thus,

        \begin{align}\label{E:lommel-nu-zero}
        G(z)& := S_{-2p-1,\,0}(L{\rm e}^{Mz}{\rm e}^{-m\pi i}) \notag\\
        &=S_{-2p-1,\,0}(L{\rm e}^{Mz})+\frac{(-1)^pK_+''}{2^{2p}(p!)^2}H_0^{(1)}(L{\rm e}^{Mz}) \\
        &\quad+\frac{(-1)^pK_-''}{2^{2p}(p!)^2}H_0^{(2)}(L{\rm e}^{Mz}) \notag
        \end{align}\smallskip

        \noindent satisfies the estimate (\ref{E:maximum-modulus}) on the
        sequence (\ref{E:sequence-z}), hence it is not subnormal for any
        integer $m$. Again the Proposition \ref{T:G(z)-growth} asserts that $S_{-2p-1,\,0}(L{\rm
        e}^{Mz})$ is not subnormal unless $K_\pm''=0$ in (\ref{E:lommel-nu-zero}), which implies that $m=0$. A contradiction to a free choice of $m$. Hence we have $K_\pm''\not=0$ and so
        $S_{-2p-1,\,0}(L{\rm e}^{Mz})$ is not subnormal.\smallskip

        The final case is treated as follows:

    \item[(iii)]  Suppose that $\nu=-k$ for a positive integer $k$. However, the Proposition \ref{T:G(z)-growth} is not applicable. Instead, we shall note from Lemma \ref{L:lommel-continuation-odd-negative-general} that if\newline $-\pi<\arg {(L{\rm e}^{Mz})}<\pi$, then for \textit{each} integer $m$,

        \begin{align}\label{E:lommel-nu-negative-integer}
       G(z)&:=S_{-k-2p-1,\,-k}(L{\rm e}^{Mz}{\rm e}^{-2m\pi i})\notag\\
        &=S_{-k-2p-1,\,-k}(L{\rm e}^{Mz})+\frac{(-1)^{k+p}}{2^{k+2p}k!p!(1+k)_p}(L{\rm e}^{Mz})^{-k} H(z),
        \end{align}\smallskip

        \noindent where the entire function $H(z)$ is defined by

        \begin{align}\label{E:H(z)-definition}
        H(z)&:=B_k(L{\rm e}^{Mz})\big[K_+''H_0^{(1)}(L{\rm e}^{Mz})+K_-''H_0^{(2)}(L{\rm e}^{Mz})] \\
        &\quad-L{\rm e}^{Mz} C_k(L{\rm e}^{Mz})\big[K_+''H_1^{(1)}(L{\rm e}^{Mz})+K_-''H_1^{(2)}(L{\rm e}^{Mz})\big]. \notag
        \end{align}\smallskip

        Now we shall obtain an estimate of the growth of $H(z)$ by following the idea of proof of the Proposition \ref{T:G(z)-growth}. We first define the polynomials \newline$D_k^\pm(\zeta):=B_k(\zeta)\pm i\zeta C_k(\zeta)$. Then we must have $D_k^\pm(\zeta) \not\equiv 0$ for $B_k(\zeta) \equiv C_k(\zeta) \equiv 0$ otherwise, which is not permitted by  {the} remark following (\ref{E:recurrence-A_n-B_n-C_n}). We further redefine the sequences (\ref{E:sequence-z}), (\ref{E:sequence-r_n}) and (\ref{E:sequence-e^z}) with the constant $C$ replaced by the polynomial $D_k^+(\zeta)$. We let $d_k^\pm=\deg D_k^\pm(\zeta)$. Here $0 \le d_k^\pm \le k+1$. Thus it follows from the asymptotic expansions (\ref{E:hankel1-asy}),  (\ref{E:hankel2-asy}) and (\ref{E:sequence-e^z}) all on the sequence (\ref{E:sequence-z}) that (\ref{E:H(z)-definition}) becomes

		\begin{align*}
        H(z_n)&=\bigg(\frac{2}{\pi L{\rm e}^{Mz_n}}\bigg)^{\frac{1}{2}}
        \Big\{D_k^+(L{\rm e}^{Mz_n})K_+''{\rm e}^{-\frac{\pi}{4}i}{\rm e}^{iL{\rm e}^{Mz_n}}\big[1+{ O((L{\rm e}^{Mz_n})^{-1})}\big]\\
        &\qquad +D_k^-(L{\rm e}^{Mz_n})K_-''{\rm e}^{\frac{\pi}{4}i}{\rm e}^{-iL{\rm e}^{Mz_n}}\big[1+{ O((L{\rm e}^{Mz_n})^{-1})}\big]\Big\}\\
        &={ \left(\frac{2}{\pi(1 \mp i)}\right)^{\frac12}R_n^{-\frac12}}
        \Big\{D^+_k(R_n(1\mp i))K_+^{\prime\prime}{\rm e}^{-i\frac{\pi}{4}}
        {\rm e}^{iR_n(1\mp i)}\big[1+{ O\big((R_n(1\mp i))^{-1}\big)}\big]\\
        &\qquad+ D^-_k(R_n(1\mp i))K_-^{\prime\prime}{\rm e}^{i\frac{\pi}{4}}
        {\rm e}^{-iR_n(1\mp i)}\big[1+{ O\big((R_n(1\mp i))^{-1}\big)}\big]\Big\}\\
		&=\left\{
			\begin{array}{ll}
				\D { \left(\frac{2}{\pi(1 - i)}\right)^{\frac12}R_n^{-\frac12}}\bigg\{
			    D^+_k(R_n(1- i))K_+^{\prime\prime} {\rm e}^{-i\frac{\pi}{4}}
				{\rm e}^{iR_n(1- i)} \\
                \qquad\times \big[1+{ O\big((R_n(1 - i))^{-1}\big)}\big]& \\
				\qquad + D^-_k(R_n(1- i))K_-^{\prime\prime} {\rm e}^{i\frac{\pi}{4}}
				{\rm e}^{-iR_n(1- i)}\big[1+{ O\big((R_n(1 - i))^{-1}\big)}\big]\bigg\},
				& \text{if\ } D^+_k\not\equiv 0;\\
                & \\
				\D { \left(\frac{2}{\pi(1 + i)}\right)^{\frac12}R_n^{-\frac12}} D^-_k(R_n(1+ i))K_-^{\prime\prime} {\rm e}^{i\frac{\pi}{4}}
				{\rm e}^{-iR_n(1+ i)} \\
                \qquad \times \big[1+{ O\big((R_n(1 + i))^{-1}\big)}\big],
				& \text{if\ } D^+_k\equiv 0.
			\end{array}\right.
        \end{align*}

    Hence
		\begin{align}\label{E:estimate-H(z)}
			|H(z_n)|& { \ge} \big({2^{1/2}}{{\pi}^{-1}}R_n^{-1}\big)^\frac12 \notag\\
			&\qquad\times
			\left\{
				\begin{array}{ll}
					|D^+_k(R_n(1- i))K_+^{\prime\prime}| {\rm e}^{R_n}
					 { -} |D^-_k(R_n(1- i))K_-^{\prime\prime}| {\rm e}^{-R_n} \\
					\qquad +o\big(R_n^{d^+_k}{\rm e}^{R_n}\big),
					& \text{if\ } D^+_k\not\equiv 0;\\
                    & \\
					|D^-_k(R_n(1+ i))K_-^{\prime\prime}| {\rm e}^{R_n}
					+o\big(R_n^{d^-_k} {\rm e}^{R_n}\big),
					& \text{if\ } D^+_k\equiv 0;\\
				\end{array}\right. \notag\\
			&=\big({2^{1/2}}{{\pi}^{-1}}R_n^{-1}\big)^\frac12 \notag\\
			&\qquad\times
			\left\{
				\begin{array}{ll}
					\D |D^+_k(R_n(1- i))K_+^{\prime\prime}| {\rm e}^{R_n}
					\Bigg(1 { -}\left|\frac{D^-_k(R_n(1- i))K_-^{\prime\prime}}
					{D^+_k(R_n(1- i))K_+^{\prime\prime}}\right| {\rm e}^{-2R_n}\Bigg)\\
					\qquad+o\big(R_n^{d^+_k}{\rm e}^{R_n}\big),
					& \text{if\ } D^+_k\not\equiv 0;\\
                    & \\
					|D^-_k(R_n(1+ i))K_-^{\prime\prime}| {\rm e}^{R_n}
					+o\big(R_n^{d^-_k}{\rm e}^{R_n}\big),
					& \text{if\ } D^+_k\equiv 0.\\
				\end{array}\right.
			\end{align}\smallskip

         Therefore it follows from the estimate (\ref{E:estimate-H(z)}) and Lemma \ref{L:lommel-modulus} that the entire function $G(z)=S_{-k-2p-1,\,-k}(L{\rm e}^{Mz}{\rm e}^{-2m\pi i})$ is not subnormal unless $K_\pm''=0$ in equation (\ref{E:lommel-nu-negative-integer}) which implies that $m=0$, a contradiction to the free choice of $m$ again. Hence we have $K_\pm''\not=0$ and so $S_{-k-2p-1,\,-k}(L{\rm e}^{Mz})$ is not subnormal. This completes the proof of the Proposition \ref{L:lommel-not-subnormal}.\smallskip

\end{enumerate}
\end{proof}
\medskip

We may now continue the proof of the Theorem 1.4 (b).\newline
To prove the result that $\mu_j$ and $\nu$ must satisfy either one of the
equations in (\ref{E:desired-equations}) when $\sigma_j\not=0$, where
$j \in \{1,\,2,\ldots,\,n\}$, we recall from the Remark \ref{R:bessel-different-branches} that

$$S_{\mu_1,\,\nu}(L{\rm e}^{Mz}),\ S_{\mu_2,\,\nu}(L{\rm
e}^{Mz}),\ldots,\ S_{\mu_n,\,\nu}(L{\rm e}^{Mz})$$\smallskip

\noindent are entire functions in the $z$-plane and that each $S_{\mu_j,\,\nu}(L{\rm e}^{Mz})$ ($j=1,\,2,\ldots,\, n$) is
\textit{independent of the branches} of $S_{\mu_j,\,\nu}(\zeta)$. This fact allows us to do the
following: Let $j$ be an element of the set $\{1,\,2,\ldots,\,n\}$
such that $\sigma_j\not=0$. For such a fixed $j$, we rewrite the solution
(\ref{E:lommel-sum}) as

\begin{align}\label{E:lommel-sum-modify}
f(z)=\sigma_j {\rm e}^{-N z} S_{\mu_j,\, \nu}(L {\rm e}^{M z}{\rm
e}^{-m\pi i})+\sum_{k=1 \atop{k\not=j}}^n \sigma_k {\rm e}^{-N
z}S_{\mu_k,\, \nu}(L {\rm e}^{M z}),
\end{align}\smallskip

\noindent where the function $S_{\mu_j,\, \nu}(\zeta)$ is in the
branch $-(m+1)\pi<\arg \zeta<-(m-1)\pi$ and the other Lommel functions $S_{\mu_1,\,\nu}(\zeta),\ldots,\,S_{\mu_{j-1},\,\nu}(\zeta), \,S_{\mu_{j+1},\,\nu}(\zeta),\ldots,\,S_{\mu_n,\,\nu}(\zeta)$ are in
the principal branch $-\pi<\arg \zeta<\pi$ and $m$ is an arbitrary but otherwise fixed non-zero integer.\smallskip

\begin{remark}\label{R:4.1}\rm
We note again that in the following discussion that we only consider the case $\mu_j-\nu=-2p_j-1$. The other case $\mu_j+\nu=-2p_j-1$ can be dealt with similarly and applying the property that each $S_{\mu_j,\,\nu}(\zeta)$ is even in $\nu$.
\end{remark}\smallskip

 If $\mu_j-\nu=-2p_j-1$ for some non-negative integer $p_j$, then it
follows from the Lemma \ref{L:lommel-modulus} and equation (\ref{E:lommel-nu-non-zero}),
(\ref{E:lommel-nu-zero}) or (\ref{E:lommel-nu-negative-integer}) in the proof of the Proposition \ref{L:lommel-not-subnormal} that $f(z)$ satisfies the estimate (\ref{E:maximum-modulus}) or (\ref{E:estimate-H(z)}) on the sequence (\ref{E:sequence-z}), thus it contradicts our assumption that $f(z)$ is subnormal. Hence $\mu_j-\nu$, and then $\mu_j+\nu$ by Remark \ref{R:4.1}, \textit{cannot} be an odd negative integer.\smallskip

Now we can apply the analytic continuation formula in the Theorem \ref{T:lommel-continuation-general} with this fixed integer $m$ to get
\smallskip

\begin{align}\label{E:lommel-continuation-generalpi}
S_{\mu_j,\, \nu}(L{\rm e}^{M z}{\rm e}^{-m\pi
i})&=K_+P_m(\cos\nu\pi,{\rm e}^{-\mu_j\pi i})H_{\nu}^{(1)}(L{\rm e}^{Mz})\\
&\quad+K_+{\rm e}^{-\nu\pi i}Q_m(\cos\nu\pi,{\rm e}^{-\mu_j\pi i})H_{\nu}^{(2)}(L{\rm e}^{Mz}) \notag\\
&\quad +(-1)^m{\rm e}^{-m\mu_j \pi i}S_{\mu_j,\, \nu}(L{\rm e}^{Mz}),\notag
\end{align}\smallskip

\noindent where $P_m(\cos\nu\pi,{\rm e}^{-\mu_j\pi i})$ and
$Q_m(\cos\nu\pi,{\rm e}^{-\mu_j\pi i})$ are polynomials as defined
in the Theorem \ref{T:lommel-continuation-general}. Then expressions
(\ref{E:lommel-sum-modify}) and (\ref{E:lommel-continuation-generalpi}) give

\begin{align}\label{E:lommel-sum-modify-further}
f(z)&=K_+\sigma_j{\rm e}^{-Nz} P_m(\cos\nu\pi,{\rm e}^{-\mu_j\pi
i})H_{\nu}^{(1)}(L{\rm
e}^{Mz})\\
&\quad+K_+\sigma_j{\rm e}^{-Nz} {\rm e}^{-\nu\pi i}Q_m(\cos\nu\pi,{\rm e}^{-\mu_j\pi
i})H_{\nu}^{(2)}(L{\rm
e}^{Mz})\notag\\
&\quad +(-1)^m\sigma_j{\rm e}^{-m\mu_j \pi i-Nz}S_{\mu_j,\, \nu}(L{\rm
e}^{Mz})+\sum_{k=1 \atop{k\not=j}}^n \sigma_k {\rm e}^{-N
z}S_{\mu_k,\, \nu}(L {\rm e}^{M z}).\notag
\end{align}\smallskip

 If either of the coefficients of $H_{\nu}^{(1)}(L{\rm e}^{Mz})$ and $H_{\nu}^{(2)}(L{\rm e}^{Mz})$ in the (\ref{E:lommel-sum-modify-further}) is non-zero, then the Proposition \ref{T:G(z)-growth} implies that the entire function $f(z)$ cannot be subnormal,  which is impossible. Thus we must have

\begin{equation}\label{E:subnormal-zero-anypi}
K_+P_m(\cos\nu\pi, {\rm e}^{-\mu_j\pi i})=0=K_+{\rm e}^{-\nu\pi i}Q_m(\cos\nu\pi,{\rm
e}^{-\mu_j\pi i}).
\end{equation}\smallskip

Now we are ready to solve the equations (\ref{E:desired-equations}), we recall again that the growth of
$S_{\mu_j,\, \nu}(L {\rm e}^{M z})$ must be independent of branches which is equivalent to equations (\ref{E:subnormal-zero-anypi}) hold for \textit{each} integer $m$. It is clear from Theorem \ref{T:lommel-continuation-general}(b) that $P_m(\cos \nu\pi,\,{\rm e}^{-\mu_j \pi i})$ and $Q_m(\cos \nu\pi,\,{\rm e}^{-\mu_j \pi i})$ \textit{cannot be both identically zero with respect to each non-zero integer} $m$, this yields from (\ref{E:subnormal-zero-anypi}) that the condition $K_+=0$ holds, \textit{i.e.}, when $\sigma_j\not=0$,

$$
\cos\bigl({\mu_j+\nu\over 2}\pi\bigr)=0\quad\textrm{or}\quad 1+{\rm
e}^{(-\mu_j+\nu)\pi i}=0,$$\smallskip

\noindent where $j \in \{1,\,2,\ldots,\, n\}$. \smallskip

 We can now complete the proof of Theorem \ref{T:Chiang-Yu}(b) by considering each
equation above individually and obtain the conclusions: either

\begin{align}\tag{1.7}
\mu_j + \nu =2p_j+1\quad \mbox{or}\quad \mu_j-\nu=2p_j+1
\end{align}\smallskip

\noindent below for non-negative integers $p_j$ when $\sigma_j\not=0$,
where $j\in \{1,\,2,\ldots,\, n\}$.\smallskip

\begin{itemize}
  \item[Case (i):] Suppose that
$\D\cos\bigl(\frac{\mu_j+\nu}{2}\pi\bigr)=0$. This equation gives

$$\mu_j+\nu=2p_j+1$$\smallskip

\noindent for some integer $p_j$. By the paragraph following Remark \ref{R:4.1}, $p_j$ must be a non-negative integer and then the Remark \ref{R:lommel-terminate} implies that the expansion of $S_{\mu_j,\,\nu}(L{\rm
e}^{Mz})$ terminates and $S_{\mu_j,\,\nu}(L{\rm
e}^{Mz})/(L{\rm e}^{Mz})^{\mu_j-1}$ becomes a polynomial in $L{\rm e}^{Mz}$ and
$1/L{\rm e}^{Mz}$, as asserted in (\ref{E:lommel-sol}).

  \item[Case (ii):] Suppose that $1+{\rm e}^{(-\mu_j+\nu)\pi i}=0$.
That is,

$$\mu_j-\nu=-2p_j-1$$\smallskip

\noindent for some integer $p_j$. Since we have shown that $\mu_j-\nu$ cannot be an odd negative integer, it follows that $p_j$ must be negative and so $\mu_j-\nu$ is an odd positive integer. Hence the Remark
\ref{R:lommel-terminate} again implies that $S_{\mu_j,\,\nu}(L{\rm
e}^{Mz})$ terminates and $S_{\mu_j,\,\nu}(L{\rm
e}^{Mz})/(L{\rm e}^{Mz})^{\mu_j-1}$ becomes a polynomial in $L{\rm e}^{Mz}$ and
$1/L{\rm e}^{Mz}$, as stated in (\ref{E:lommel-sol}).
\end{itemize}\smallskip

We recall that $j$ is an \textit{arbitrary} element in the set $\{1,\,2,\ldots,\,n\}$ such that $\sigma_j \not=0$. Thus the above argument is valid for each such $j$ and hence we have the necessary part of the Theorem.
\medskip

Conversely, suppose $A=B=0$ for $f(z)$ in (\ref{E:Chiang-Yu-soln-1})
and either $\mu_j + \nu=2p_j+1$ or $\mu_j-\nu=2p_j+1$, $p_j$ is a non-negative integer,
where $j\in \{1,\,2,\ldots,\, n\}$ with $\sigma_j\not=0$. Then
clearly, (according to the Remark \ref{R:lommel-terminate} that) each
$S_{\mu_j,\,\nu}(L{\rm e}^{Mz})/(L{\rm e}^{Mz})^{\mu_j-1}$ is a polynomial in $L{\rm e}^{Mz}$
and/or $1/L{\rm e}^{Mz}$. Hence $f(z)$ is clearly subnormal. This
proves the converse part and so completes the proof of Theorem \ref{T:Chiang-Yu}(b).\hfill $\square$\medskip

\section{\bf Proof of Theorem {\rm \ref{T:lommel}} and a consequence}

The proof is a direct consequence of the proof to the Theorem \ref{T:Chiang-Yu} given in \S\ref{P:Chiang-Yu}. In fact, the Proposition \ref{T:G(z)-growth} asserts that

\begin{equation}\label{E:G(z)-2}
G(z)=\sum_{j=1}^n \sigma_jS_{\mu_j,\, \nu}(L{\rm e}^{Mz})
\end{equation}\smallskip

\noindent whenever $G(z)$ is of finite order of growth. The argument in remaining proof of the Theorem \ref{T:Chiang-Yu} that $\Re(\mu_j)$ are distinct certainly applies when we have only a \textit{single} $S_{\mu,\,\nu}(L{\rm e}^{Mz})$. Thus we must have either

\[
	\mu - \nu=2p+1\quad\textrm{or}\quad \mu+\nu=2p+1
\]\smallskip

\noindent for a non-negative integer $p$. Conversely, we suppose that $\mu-\nu=2p+1$. Then the Remark \ref{R:lommel-terminate} gives that

$$S_{\nu+2p+1,\, \nu}(L{\rm e}^{M z})/(L{\rm e}^{Mz})^{\nu+2p}$$ \smallskip

\noindent is the composition of a polynomial and the exponential, and hence it is of finite order of growth. Since the entire function $(L{\rm e}^{Mz})^{\nu+2p}$ is certainly of finite order of growth, we have the result that the function $S_{\nu+2p+1,\, \nu}(L{\rm e}^{M z})$ is also of finite order of growth. The case when $\mu+\nu=2p+1$ now follows easily from Remark \ref{R:4.1}. This completes the proof of Theorem \ref{T:lommel}. \hfill $\square$\smallskip

 Moreover, it follows from (\ref{E:struve}) {and (\ref{E:struve-analytic-continuation-formula})} that { for every integer $m$},

\begin{align*}
{{\bf H}_\nu({\rm e}^z {\rm e}^{-m\pi i})}&={ (-1)^m{\rm e}^{-m\nu\pi i}} {\bf H}_\nu({\rm e}^z)\\
&={ (-1)^m{\rm e}^{-m\nu\pi i}}\bigg[Y_\nu({\rm e}^z)+\frac{2^{1-\nu}}{\sqrt{\pi}\Gamma(\nu+\frac{1}{2})}
S_{\nu,\,\nu}({\rm e}^z)\bigg]\\
&={ (-1)^m{\rm e}^{-m\nu\pi i}}\bigg\{\frac{1}{2i}\big[H_\nu^{(1)}({\rm e}^z)-H_\nu^{(2)}({\rm e}^z)\big]+\frac{2^{1-\nu}}{\sqrt{\pi}\Gamma(\nu+\frac{1}{2})}
S_{\nu,\,\nu}({\rm e}^z)\bigg\}.
\end{align*}\smallskip

\noindent Then the Proposition \ref{T:G(z)-growth} immediately implies the new result:

\begin{corollary}\label{C:Struve} Let $\nu$ be an arbitrary complex number. The composition of the Struve function (irrespective of branches) and the exponential function (which is an entire function)
${\bf H}_\nu({\rm e}^z)$ is of infinite order of growth. In particular, it is not subnormal.
\end{corollary}

Let $f(z)$ be an entire function of order
$\sigma$, where $0<\sigma<1$. Then it is easy to check by the
definitions that the entire function $g(z)={\rm e}^{f(z)}$ is subnormal and has an infinite order of growth. Then Corollary \ref{C:Struve} shows that the ${\bf H}_\nu({\rm e}^z)$ grows faster than subnormal solutions.
\medskip

\section{\bf Quantization-type results and examples}
\label{S:1/16-results}
Ismail and one of the authors strengthened  \cite{CI02(2)} (announced in \cite{CI02(1)}) earlier results of Bank, Laine and Langley \cite{BL83}, \cite{BLL86} (see also \cite{Shimomura2002}) that an entire solution $f$ of either the equation

\begin{equation}
\label{E:Chiang-Ismail-1}
	y^{\prime\prime} +{\rm e}^z y=Ky
\end{equation}\smallskip

\noindent or

\begin{equation}
\label{E:Chiang-Ismail-2}
	y^{\prime\prime}+\Big(-\frac14 {\rm e}^{-2z}+\frac12{\rm e}^{-z}\Big)y=Ky
\end{equation}\smallskip

\noindent can be solved in terms of Bessel functions and Coulomb Wave functions respectively. Besides, we identify two classes of classical orthogonal polynomials (Bessel and generalized Bessel polynomials respectively) in the explicit representation of solutions under the assumption $\D\lambda(f)=\lim_{r\to+\infty}\log n(r)/\log r<+\infty$ (\textit{boundary condition}). This also results in a complete determination of the eigenvalues and eigenfunctions of the equations. This is known as the semi-classical quantization in quatum mechanics. Both equations have important physical applications. For example, the Eqn. (\ref{E:Chiang-Ismail-1})  is derived as a reduction of a non-linear Schr\"odinger equation in a recent study of Benjamin-Feir instability pheonmena in deep water in \cite{Segur2005}, while the second Eqn. (\ref{E:Chiang-Ismail-2}) is a standard classical diatomic model in quantum mechanics introduced by P. M. Morse in 1929 \cite{Morse1929}\begin{footnote} { See \cite[pp. 1-4]{slater} for a historical background of the Morse potential.}\end{footnote}. However, the equation (\ref{E:Chiang-Ismail-2}) also appears as a basic model in the recent $\mathcal{PT}-$symmetric quantum mechanics research \cite{Znojil1999} (see also \cite{Bender2007}).
\smallskip

We now consider special cases of the Theorem \ref{T:Chiang-Yu} so that the equations (\ref{E:Chiang-Yu}) exhibits
a kind of semi-classical quantization phenomenon that usually only applies to homogeneous equations. In particular, these equations admits classical polynomials solutions (Neumann's polynomials, Gegenbauer's polynomials, Schl\"{a}fli's polynomials and Struve's functions) that are related to special functions when the first derivative term in
(\ref{E:Chiang-Yu}) is zero.\smallskip

Suppose that $L=2,\, M=\frac{1}{2},\,N=0$ and $n=1$ in Theorem
{\rm \ref{T:Chiang-Yu}} so that the differential equation {\rm
(\ref{E:Chiang-Yu})} becomes

\begin{equation}\label{E:reduced-Chiang-Yu}
f''+({\rm e}^z-K)f=\sigma 2^{\mu-1}{\rm e}^{\frac{1}{2}(\mu+1)z},
\end{equation}\smallskip

\noindent where $\D K=\frac{\nu^2}{4}.$\smallskip

\begin{theorem}
\label{T:quantizations} Then, in each of the
cases below, we have the necessary and sufficient condition on $K$ that depends on the non-negative integer $p$
so that the equation {\rm (\ref{E:reduced-Chiang-Yu})} admits a
subnormal solution. Furthermore, the forms of the subnormal
solutions are given explicitly in Table \ref{mytable}:
\newline

\begin{table}[h]\vspace*{-3ex}\extrarowheight=7pt
\caption[]{Special cases of {\rm (\ref{E:reduced-Chiang-Yu})}.} \label{mytable}
\vspace{0cm}
\begin{tabular}{cccc}
   \hline
   & Cases &  Corresponding $K$ & Subnormal solutions \\ \hline
  (1)& $\mu=1$ & $p^2$ & $2\sigma {\rm e}^{\frac{z}{2}}O_{2p}(2{\rm e}^{\frac{z}{2}})$ \\
  (2)& $\mu=0$ & $\D\frac{(2p+1)^2}{4}$  & $\D \frac{2\sigma}{2p+1}{\rm e}^{\frac{z}{2}}O_{2p+1}(2{\rm e}^{\frac{z}{2}})$\\
  (3)& $\mu=-1$ & $(p+1)^2$   & $\D \frac{\sigma}{4(p+1)}S_{2p+2}(2{\rm e}^{\frac{z}{2}})$ \\
  (4)& $\mu=\nu$ &  $\D \frac{(2p+1)^2}{16}$ & $\sigma 2^{p-\frac{1}{2}}\sqrt{\pi}p!\Big[{\bf H}_{p+\frac{1}{2}}(2{\rm e}^{\frac{z}{2}})-Y_{p+\frac{1}{2}}(2{\rm e}^{\frac{z}{2}})\Big]$ \vspace{3pt}\\ \hline
\end{tabular}
\end{table}
\end{theorem}\smallskip

Here $O_{2p}(\zeta)$ and $O_{2p+1}(\zeta)$ are the Neumann
polynomials of degrees $2p$ and $2p+1$ respectively; $S_p(\zeta)$
is the Schl\"{a}fli polynomial and ${\bf
H}_{p+\frac{1}{2}}(\zeta)$ is the Struve function. Their
background information will be given in Appendix \ref{Adx:Special}.\smallskip

\begin{proof}[Proof of Theorem \ref{T:quantizations}] We only prove the first case and the other cases can be
dealt with similarly. We note that $O_{2p}(\zeta)$ is actually a
polynomial in $1/\zeta$ and that the degree of each individual
term is odd, so we see that $2\sigma {\rm
e}^{\frac{z}{2}}O_{2p}(2{\rm e}^{\frac{z}{2}})$ indeed has the
form ${\rm e}^{dz}S({\rm e}^z)$ where $S(\zeta)$ is a polynomial
as asserted in (\ref{E:Gundersen-Steinbart-soln}). According to
Theorem \ref{T:Chiang-Yu}, $f(z)$ is subnormal if and only if $\mu
\pm \nu=2p+1$, \textit{i.e.}, $\nu^2=4p^2$, where $p$ is a non-negative integer and thus

$$K=\frac{\nu^2}{4}=p^2.$$\smallskip

\noindent It follows from (\ref{E:neumann}) that

\begin{equation*}
f(z)=\sigma S_{1,\, 2p}(2{\rm e}^{\frac{z}{2}})=2\sigma {\rm
e}^{\frac{z}{2}}O_{2p}(2{\rm e}^{\frac{z}{2}}).
\end{equation*}\smallskip
\end{proof}

\begin{remark}\rm In fact, we also have

$$f(z)=\sigma S_{1,\, 2p}(2{\rm e}^{\frac{z}{2}})=\sigma
{\rm e}^{\frac{z}{2}}A_{2p,\,0}(2{\rm e}^{\frac{z}{2}}),$$\smallskip

\noindent where the $A_{2p,\,0}(\zeta)$ is Gegenbauer's polynomial
which will be discussed in \S \ref{Adx:Special}.2.
\end{remark}

\begin{example}[Even Neumann's polynomial]\label{Ex:example-even-neumann}
 \rm Take $p=1$ in Theorem \ref{T:quantizations}, Case (1), then we have

$$f(z)=\sigma S_{1,\, 2}(2{\rm e}^{\frac{z}{2}})=2\sigma {\rm e}^{\frac{z}{2}}O_{2}(2{\rm e}^{\frac{z}{2}})=\sigma +\sigma
{\rm e}^{-z}$$\smallskip

\noindent which is a subnormal solution of the equation

$$f''+({\rm e}^z-1)f=\sigma {\rm e}^z.$$
\end{example}

\begin{example}[Odd Neumann's polynomial]\label{Ex:example-odd-neumann}\rm Take $p=1$ in Theorem \ref{T:quantizations}, Case (2), then we have

$$f(z)=\sigma S_{0,\, 3}(2{\rm e}^{\frac{z}{2}})= {2\sigma\over 3}
{\rm e}^{\frac{z}{2}}O_3(2{\rm
e}^{\frac{z}{2}})=\frac{\sigma}{2}{\rm e}^{-\frac{z}{2}}+\sigma
{\rm e}^{-\frac{3z}{2}}$$\smallskip

\noindent which is a subnormal solution of the equation

$$f''+\bigg({\rm e}^z-\frac{9}{4}\bigg)f=\frac{\sigma}{2} {\rm e}^{\frac{z}{2}}.$$
\end{example}

\begin{example}[Schl\"{a}fi's polynomial] \rm Take $p=1$ in Theorem \ref{T:quantizations}, Case (3), then we have

$$f(z)=\sigma S_{-1,\, 4}(2{\rm e}^{\frac{z}{2}})=
{\sigma \over 8} S_{4}(2{\rm
e}^{\frac{z}{2}})=\frac{\sigma}{4}{\rm
e}^{-z}+\frac{3\sigma}{4}{\rm e}^{-2z}$$\smallskip

\noindent which is a subnormal solution of the equation

$$f''+({\rm e}^z-4)f=\frac{\sigma}{4}.$$
\end{example}

\begin{example}[Struve's function]\label{Ex:example-struve}\rm Take $p=1$ in Theorem \ref{T:quantizations}, Case (4), then we have

\begin{align*}
f(z)=\sigma S_{\frac{3}{2},\frac{3}{2}}(2{\rm
e}^{\frac{z}{2}})=\sigma \sqrt{2}\sqrt{\pi} \bigg[{\bf
H}_{\frac{3}{2}}(2{\rm e}^{\frac{z}{2}})-Y_{\frac{3}{2}}(2{\rm
e}^{\frac{z}{2}})\bigg]=\sigma \sqrt{2} {\rm e}^{\frac{z}{4}}
\bigg(1+\frac{1}{2{\rm e}^z}\bigg)
\end{align*}\smallskip

\noindent which is a subnormal solution of the equation

$$f''+\bigg({\rm e}^z-\frac{9}{16}\bigg)f=\sigma \sqrt{2} {\rm e}^{\frac{5z}{4}}.$$
\end{example}

\begin{remark} \rm We note that the differential equations in Examples
\ref{Ex:example-odd-neumann} and \ref{Ex:example-struve} are not in the form (\ref{E:Gundersen-Steinbart})
and therefore their corresponding solutions are not in the form
(\ref{E:Gundersen-Steinbart-soln}). These give explicit examples
to Remark \ref{R:relation-theoremA-ChiangYu}.
\end{remark}\medskip

\section{\bf Concluding remarks}

Semi-explicit representations of special entire solutions for homogeneous second order linear periodic differential equations were obtained by complex analysts using tools from Nevanlinna theory starting from the 1970s (see for examples \cite{Fr61}, \cite{Wi67}, \cite{BL83}). The main idea has an origin in the Picard$\backslash$Borel exceptional value for entire functions and the theories develop quite independent from the classical approaches to Hill's equations and Mathieu equations (\cite{WW27}). Many of these equations are related to models in quantum mechanics (\cite{Morse1929}) and it is found that the determination of their (discrete) eigenvalues (\textit{i.e.,} their spectrum) can be described by using exponent of convergence of the zeros of solutions, a quantity that is related to Borel exceptional value in classical function theory (\cite{BL83}, \cite{CI02(2)}). Moreover, one of the authors and Ismail identified an important class of orthogonal polynomials: generalized Bessel polynomials appear in the representation of the special solutions to (\ref{E:Chiang-Ismail-1}) and (\ref{E:Chiang-Ismail-2}) in \cite{CI02(2)} along with the determination of their eigenvalues. This paper extends, on the one hand, that one can use the Lommel functions and related polynomials (to Bessel functions) to describe the subnormal (special) solutions of (\ref{E:Chiang-Yu-special}) first considered by Gundersen and Steinbart \cite{GS94}, and on the other hand, shows that a kind of semi-classical quantization for non-homogeneous equations also exists for (\ref{E:Chiang-Yu-special}), (\ref{E:Chiang-Yu}) and the equations in the Theorem \ref{T:quantizations}. In addition, we obtain a number of new analytic continuation formulae for the Lommel functions, and a new property for the Lommel functions (Theorem \ref{T:lommel}). Although the Lommel functions have numerous physical applications as mentioned in the Introduction, to the best of the authors' knowledge, only few papers have been written to investigate their mathematical properties in the past decades. See for examples \cite{Dingle1}, \cite{Dingle2}, \cite{Steinig1972}, \cite{Pidduck1946}, \cite{Rollinger1964} and \cite{Goard2003}, and the references therein.\smallskip

Although we generally do not have a simple quantum mechanical interpretation for non-homogeneous equations like the equation (\ref{E:Chiang-Yu-special}),  its homogeneous counterpart  and the Lommel functions themselves have numerous applications in various branches of physical applications. So it is hoped the results in this paper will be of interest for others in due course.

\medskip

\appendix

\section{\bf Preliminaries on Bessel functions}\label{Adx:Bessel}
 Let $m$ be an integer. We record here the following analytic continuation formulae for the Bessel functions
\cite{Wa44}, \S3.62:

\begin{align}
J_\nu(\zeta {\rm e}^{m\pi i})&= {\rm e}^{m\nu \pi i} J_\nu(\zeta),
\label{E:bessel-first-continuation}\\
Y_\nu(\zeta {\rm e}^{m\pi i})&= {\rm e}^{-m\nu \pi i} Y_\nu
(\zeta)+ 2i \sin(m\nu\pi) \cot(\nu\pi)
J_\nu(\zeta).\label{E:bessel-second-continuation}
\end{align}\smallskip

\noindent We recall the {\it Bessel functions of the third kind of order} $\nu$ \cite{Wa44}, \S3.6 are given by

\begin{equation}\label{E:hankels-definition}
H_\nu^{(1)}(\zeta)=J_\nu(\zeta)+iY_\nu(\zeta), \quad
H_\nu^{(2)}(\zeta)=J_\nu(\zeta)-iY_\nu(\zeta).
\end{equation}\smallskip

\noindent They are also called the \textit{Hankel functions of order $\nu$ of the first and second kinds}.\smallskip

The asymptotic expansions of $H_\nu^{(1)}(\zeta)$ and $H_\nu^{(2)}(\zeta)$ are also recorded as follows:

\begin{equation}
\label{E:hankel1-asy} \left(\frac{\pi
\zeta}{2}\right)^{\frac{1}{2}}H_\nu^{(1)}(\zeta)= {\rm
e}^{i\left(\zeta-\frac{1}{2}\nu\pi-\frac{1}{4}\pi\right)}
\left[\sum_{k=0}^{p-1}\frac{\left(\frac{1}{2}-\nu\right)_k\left(\frac{1}{2}
+\nu\right)_k} {k!(2i\zeta)^k}+R_p^{(1)}(\zeta)\right]
\end{equation}\smallskip

\noindent where $R_p^{(1)}(\zeta)=O(\zeta^{-p})$ in $-\pi < \arg \zeta < 2\pi$;

\begin{equation}
\label{E:hankel2-asy} \left(\frac{\pi
\zeta}{2}\right)^{\frac{1}{2}}H_\nu^{(2)}(\zeta)= {\rm
e}^{-i\left(\zeta-\frac{1}{2}\nu\pi-\frac{1}{4}\pi\right)}
\left[\sum_{k=0}^{p-1}\frac{\left(\frac{1}{2}-\nu\right)_k\left(\frac{1}{2}
+\nu\right)_k} {k!(-2i\zeta)^k}+R_p^{(2)}(\zeta)\right]
\end{equation} \smallskip

\noindent where $R_p^{(2)}(\zeta)=O(\zeta^{-p})$ in $-2\pi < \arg \zeta < \pi$. See \cite{Wa44}, \S7.2. As a result, we see that the asymptotic expansions (\ref{E:hankel1-asy}) and (\ref{E:hankel2-asy}) are valid \textit{simultaneously} in the range $-\pi<\arg \zeta<\pi$.\smallskip

Now we deduce from (\ref{E:bessel-first-continuation}) and (\ref{E:bessel-second-continuation}) the analytic continuation formulae for $H_\nu^{(1)}(\zeta)$ and $H_\nu^{(2)}(\zeta)$
\cite{Wa44}, \S3.62:

\begin{align}
H_\nu^{(1)}(\zeta {\rm e}^{m\pi i})&=\frac{\sin(1-m)\nu\pi}{\sin
\nu\pi} H_\nu^{(1)}(\zeta)-\frac{{\rm e}^{-\nu\pi i}\sin
m\nu\pi}{\sin\nu\pi}H_\nu^{(2)}(\zeta),\label{E:bessel-third-1-continuation}\\
H_\nu^{(2)}(\zeta {\rm e}^{m\pi i})&=\frac{{\rm e}^{\nu\pi i}\sin
m\nu\pi}{\sin\nu\pi}
H_\nu^{(1)}(\zeta)+\frac{\sin(m+1)\nu\pi}{\sin\nu\pi}H_\nu^{(2)}(\zeta).
\label{E:bessel-third-2-continuation}
\end{align} \smallskip

\noindent We note that the right hand sides of (\ref{E:bessel-third-1-continuation}) and
(\ref{E:bessel-third-2-continuation}) are the principal branch of the Hankel functions.\smallskip

Finally, we record the following derivative formulae for the Hankel functions \cite{Wa44}, \S3.6:

\begin{equation}\label{E:hankel-derivatives}
\frac{{\rm d}}{{\rm d\zeta}}H_0^{(1)}(\zeta)=-H_1^{(1)}(\zeta),\quad \frac{{\rm d}}{{\rm d\zeta}}H_0^{(2)}(\zeta)=-H_1^{(2)}(\zeta).
\end{equation}

\section{\bf Special polynomials and functions}\label{Adx:Special}
\subsection{Neumann's polynomials} The Neumann polynomials $O_n(\zeta)$ are defined as the
coefficients in the expansion of $\D \frac{1}{\zeta-z}$ in terms of Bessel functions $J_j(z)$ with $j \ge 0$, see \cite{Wa44}, chap. IX and the connection between the Lommel function $S_{\mu,\,\nu}(\zeta)$ and Neumann's polynomials is given by \cite{Wa44}, \S9.1 and \S10.74:


\begin{align}\label{E:neumann}
O_n(\zeta)&=\frac{1}{4}\sum_{m=0}^n \frac{n\Gamma(\frac{1}{2}n+\frac{1}{2}m)\cos^2\frac{1}{2}(m \pm n)\pi}{\Gamma(\frac{1}{2}n-\frac{1}{2}m+1)}\bigg(\frac{\zeta}{2}\bigg)^{-m-1}\\
&=\left\{
             \begin{array}{ll}
               \D \frac{S_{1,\,2p}(\zeta)}{\zeta}, \hspace{2.6cm}\hbox{if $n=2p$;} \\
               \D \frac{(2p+1)S_{0,\,2p+1}(\zeta)}{\zeta}, \hspace{1cm}\hbox{if $n=2p+1$.} \\
             \end{array}
           \right. \notag
\end{align}

\subsection{Gegenbauer's polynomials} The Gegenbauer polynomials $A_{n,\,\nu}(\zeta)$ appear in the problem of expansion of $\D \frac{\zeta^\nu}{t-\zeta}$. The formulae connecting Lommel's function $S_{\mu,\,\nu}(\zeta)$ and Gegenbauer's polynomial $A_{n,\,\nu}(\zeta)$ are given by

\begin{align*}
A_{n,\,\nu}(\zeta)&=\frac{2^{\nu+n}}{\zeta^{n+1}}\sum_{m=0}^{\le \frac{1}{2}n}\frac{\Gamma(\nu+n-m)}{m!}\bigg(\frac{\zeta}{2}\bigg)^{2m}\\
&=\left\{
             \begin{array}{l}
             \D\frac{2^\nu\Gamma(\nu+p)}{p!}\cdot\frac{\nu+2p}{\zeta^{1-\nu}}S_{1-\nu,\,\nu+2p}(\zeta), \hspace{2.1cm} \hbox{if $n=2p$;} \\
             \D\frac{2^{\nu+1}\Gamma(\nu+p+1)}{p!}\cdot\frac{\nu+2p+1}{\zeta^{1-\nu}}S_{-\nu,\,\nu+2p+1}(\zeta),\quad
             \hbox{if $n=2p+1$.}
             \end{array}
           \right.
\end{align*}\smallskip

\noindent See \cite{Wa44}, \S9.2 and \S10.74 for details.

\subsection{Schl\"afli's polynomials} The polynomial is defined by
$S_0(\zeta)=0$ and
\begin{equation*} {1\over 2}pS_p(\zeta)=\zeta O_p(\zeta)-\cos^2\biggl({p\pi\over 2}\biggr)
\end{equation*}\smallskip

\noindent for each positive integer $p \ge 1$. The polynomials are mainly introduced because they have greater simplicity over Neumann's polynomials \cite{Wa44}, \S9.3. The explicit formulae are
given by

\begin{align}\label{E:schlafli}
S_n(\zeta)=\left\{
             \begin{array}{ll}
               \D \sum_{m=1}^p {(p+m-1)! \over (p-m)!}\biggl(\frac{\zeta}{2}\biggr)^{-2m}, & \hbox{if $n=2p$;} \\
               \D \sum_{m=1}^p {(p+m)! \over (p-m)!}\biggl(\frac{\zeta}{2}\biggr)^{-2m-1}, & \hbox{if $n=2p+1$.}
             \end{array}
           \right.
\end{align}

\subsection{Lommel's functions} We briefly discuss the nature of the Lommel function $S_{\mu,\,\nu}(\zeta)$ that we encountered. We start with an equation of more general form

\begin{equation*}
\label{E:} w^{\prime\prime}+\left[\sum_{j=0}^{+\infty} {f_j\over
x^j}\right]w^\prime+\left[\sum_{j=0}^{+\infty} {g_j\over
x^j}\right]w=z^\alpha \left[\sum_{j=0}^{+\infty} {p_j\over
x^j}\right],
\end{equation*}\smallskip

\noindent in $|x|>a$ for some $a>0$. The solution to this equation
is given by

\[
w(x)=A w_1(x)+B w_2(x)+W(x),
\]\smallskip

\noindent where $A,B \in \mathbf{C}$. The asymptotic expansions of the complementary functions $w_1$ and $w_2$ are considered in \cite{Ol97}, chap. 7. It can be shown (\cite{Ol97}, chap. 7) that there exists a particular integral $W$ so that its asymptotic expansion is given by

\begin{equation}
\label{E:slow-expansion}
W_n(x)=x^\alpha\left[\sum_{j=0}^{n-1}{a_j\over
x^j}+O\biggl({1\over x^n}\biggr)\right],
\end{equation}\smallskip

\noindent as $x\to \infty$ in an unbounded region that depends on the coefficients $f_0$ and $g_0$. We note that there is no claim that such a particular integral is unique. We refer the reader to \cite{Ol97} for the details.\smallskip

The Lommel function $s_{\mu,\,\nu}(\zeta)$ given by (\ref{E:small-lommel}) is a ${}_1 F_2$ function. The
$S_{\mu,\,\nu}(\zeta)$ defined by (\ref{E:lommel-lommel-1}) and (\ref{E:lommel-lommel-2}) matches the $W(x)$ mentioned above, as clearly indicated by its asymptotic expansion in (\ref{E:lommel-asymptotic-expansion}). Thus the Lommel functions contribute to the subnormal solutions that concerns us in this paper.

\subsection{Struve's functions} If we take $n=1$, we have $\sigma_1=[2^{\nu-1}\sqrt{\pi}\Gamma\left(\nu+\frac{1}{2}\right)]^{-1}$ and $\mu=\nu$ in equation (\ref{E:generalized-Lommel}), then $\sigma_1 s_{\nu,\, \nu}(\zeta)$ is the particular solution of it.
Furthermore, it is known that the Struve function ${\bf H}_\nu(\zeta)$ of order $\nu$ \cite{Bateman}, pp. 37-39 which is defined by

$${\bf H}_\nu(\zeta):=\bigg(\frac{1}{2}\zeta\bigg)^{\nu+1}\sum_{m=0}^{+\infty}
\frac{(-1)^m\big(\frac{1}{4}\zeta^2\big)^m}{\Gamma(m+\frac{3}{2})\Gamma(\nu+m+\frac{3}{2})},$$\smallskip

\noindent also satisfies the same differential equation

\begin{equation}\label{E:struve-eqn}
\zeta^2 y''(\zeta)+\zeta y'(\zeta)+(\zeta^2-\nu^2)y(\zeta)
=\sigma_1 \zeta^{\nu+1}.
\end{equation}\smallskip

Following Olver's notation \cite{Ol97}, chap. 7, we call

\begin{equation*}
\textbf{K}_\nu(\zeta)=\textbf{H}_\nu(\zeta)-Y_\nu(\zeta)
\end{equation*}\smallskip

\noindent and it follows from (\ref{E:lommel-lommel-1}) that

\begin{align}\label{E:struve}
\textbf{H}_\nu(\zeta)={2^{1-\nu}
\over{\sqrt{\pi}\Gamma(\nu+\frac{1}{2})}}s_{\nu,\,\nu}(\zeta) =
Y_\nu(\zeta)+{2^{1-\nu} \over
\sqrt{\pi}\Gamma(\nu+\frac{1}{2})}S_{\nu,\,\nu}(\zeta).
\end{align}\smallskip

{By (\ref{E:struve}), the continuation formula for the Struve function can be determined immediately. In fact, when $\nu \neq -n-\frac{1}{2}$ for any non-negative integer $n$, then} it is noted easily from (\ref{E:small-lommel-analytic-continuation}) and {the first half of} (\ref{E:struve}) that 

\begin{equation}\label{E:struve-analytic-continuation-formula}
\textbf{H}_\nu(\zeta {\rm e}^{-m\pi i})=(-1)^m{\rm e}^{-m\nu\pi i}\textbf{H}_\nu(\zeta),
\end{equation}\smallskip

\noindent where $m$ is any integer. {When $\nu=-n-\frac{1}{2}$ for a non-negative integer $n$, then $2\nu=-2n-1$ is an odd negative integer and so $s_{\nu,\,\nu}(\zeta)$ is undefined in this case but the second half of (\ref{E:struve}) gives $\textbf{H}_{-n-\frac{1}{2}}(\zeta)=Y_{{-n-\frac{1}{2}}}(\zeta)=(-1)^nJ_{n+\frac{1}{2}}(\zeta)$. Then it follows from (\ref{E:bessel-first-continuation}) that the analytic continuation formula (\ref{E:struve-analytic-continuation-formula}) is also valid at $\nu=-n-\frac{1}{2}$. (We shall note that (\ref{E:struve-analytic-continuation-formula}) was already given in \cite{Wa44}, \S 10.41 (5).)}\smallskip

It is shown \cite{Ol97}, chap. 7 that the function $\textbf{K}_\nu(\zeta)$ defined above is the \textit{unique} particular integral of the Struve equation that satisfies (\ref{E:slow-expansion}). The expansion terminates if and only if $\nu$ is half of an odd positive integer. Example \ref{Ex:example-struve} corresponds to this situation. So the general solution for the equation (\ref{E:struve-eqn}) can be written in terms of Bessel functions and $\textbf{K}_\nu(\zeta)$.\smallskip

We remark that we can prove a special case directly on (\ref{E:struve-eqn}) by appealing to the analytic continuation
formula \cite{Ol97}, chap. 7, Ex. 15.4

\[
\textbf{K}_\nu(\zeta {\rm e}^{-\pi i})=-{\rm e}^{-\nu\pi
i}\textbf{K}_\nu(\zeta)+2i\cos(\nu\pi) H_\nu^{(1)}(\zeta),
\]\smallskip

\noindent which is a special case of Lemma \ref{L:lommel-continuation-hankel}.
\medskip

\noindent\textbf{Acknowledgements} We would like to thank Mourad E.
H. Ismail for useful conversation concerning Lemma \ref{E:lommel-independence-1}. We would also like to thank the referee for his/her careful reading of our manuscript.\smallskip

    \begin{note} {\rm When the authors revised the paper, we realized that the following simple relation holds for $P_m(\cos \nu\pi,\,{\rm e}^{-\mu\pi i})$ and $Q_m(\cos \nu\pi,\,{\rm e}^{-\mu\pi i})$ in the Theorem 3.4:

    \begin{equation*}
    Q_m(\cos \nu\pi,\,{\rm e}^{-\mu\pi i})=P_{m-1}(\cos \nu\pi,\,{\rm e}^{-\mu\pi i}),
    \end{equation*}\smallskip

    \noindent for a non-zero integer $m$. This simplifies the statement of the Theorem.}
    \end{note}

\end{document}